\documentclass[11pt,a4paper]{amsart}
\usepackage[T1]{fontenc}
\usepackage{lmodern,amssymb,booktabs}
\usepackage[a4paper,left=34mm,right=34mm,top=29mm,bottom=29mm]{geometry}
\usepackage{needspace}
\usepackage{hyperref}
\hypersetup{hidelinks,pdfauthor={Weidong Gao, Alfred Geroldinger, Guoqing Wang and Xiangneng Zeng},
 pdftitle={A disproof of the van Emde Boas--Kruyswijk conjecture on the Davenport constant},
 pdfkeywords={Davenport constant, zero-sumfree sequence, homocyclic group}}
\newcommand{\DD}{\mathsf D}
\renewcommand{\keywordsname}{Keywords}
\makeatletter
\def\@setkeywords{{\itshape\keywordsname:}\enspace\@keywords\@addpunct.}
\makeatother
\newcommand{\paperlink}[2]{\href{#1}{\textit{#2}}}
\newtheorem{conjecture}{Conjecture}
\newtheorem{theorem}{Theorem}[section]
\newtheorem{lemma}[theorem]{Lemma}
\newtheorem*{auxiliarylemma}{{Lemma}}
\theoremstyle{remark}
\newtheorem*{remark}{Remark}

\title[The van Emde Boas--Kruyswijk Conjecture]{A Disproof of the van Emde Boas--Kruyswijk Conjecture on the Davenport constant}
\author{Weidong Gao}
\address{Center for Applied Mathematics and KL-AAGDM \\ Tianjin University \\ Tianjin 300072, P.R. China}
\email{wdgao@nankai.edu.cn}
\author{Alfred Geroldinger}
\address{Department of Mathematics and Scientific Computing \\ University of Graz \\ Graz, Austria}
\email{alfred.geroldinger@uni-graz.at}
\author{Guoqing Wang}
\address{School of Mathematical Sciences \\ Tiangong University \\ Tianjin 300387, P.R. China}
\email{gqwang1979@aliyun.com}
\author{Xiangneng Zeng}
\address{Sino-French Institute of Nuclear Engineering and Technology,
Sun Yat-Sen University, Guangzhou 510275, P.R. China}
\email{zengxneng@mail.sysu.edu.cn}
\date{}

\newcommand{\D}{\DD}
\theoremstyle{plain}
\newtheorem{proposition}[theorem]{Proposition}
\theoremstyle{remark}
\begin{document}
\begin{abstract}
The Davenport constant $\DD(G)$ of a finite abelian group $G$ is the smallest
positive integer $\ell$ such that every sequence $S$ with terms from $G$ and with length  $|S| \ge \ell$
has a nontrivial zero-sum subsequence.  A long-standing
conjecture, first proposed by van Emde Boas and Kruyswijk in 1969, states that
$\DD(C_n^r)=r(n-1)+1$ for all positive integers $n$ and $r$. The equality holds true if $n$ is a prime power or if $r \le 2$. We disprove this conjecture for all
exponents $n>1$ that are not prime powers.
On the other hand, for every fixed rank $r\ge3$, we prove that $\DD(C_n^r)=r(n-1)+1$ for all sufficiently large integers $n$.
\end{abstract}
\maketitle

\section{Introduction}
Let $G = C_{n_1} \oplus \ldots \oplus C_{n_r}$ with $1 < n_1 \mid \ldots \mid n_r$ be a finite abelian group.  A straightforward example shows that
\begin{equation} \label{lower-bound}
\mathsf D^* (G) :=  1 + \sum_{i=1}^r (n_i-1) \le \mathsf D (G) \,.
\end{equation}
Already in the 1960s, it was proved that equality holds for $p$-groups and for groups of rank at most two. This was achieved independently by Kruyswijk and Olson \cite{olsonI,olsonII}. Also, already in the 1960s, a first counterexample was found showing that strict inequality may hold.
Since the Davenport constant is closely connected to a variety of problems in graph theory, finite geometry, combinatorics, and factorization theory, it received wide attention in the last decades. In parallel with the study of the exact value of the Davenport constant, the structure of minimal zero-sum sequences of length $\DD (G)$ was investigated. There is a wealth of literature on this topic (see \cite{Ge-Gr-Zh26a} for an overview).

In 1969, van Emde Boas and Kruyswijk \cite[\S12, p.~31]{vk} proposed
the following conjecture.

\begin{conjecture}\label{conj:homocyclic}
For all positive integers $n$ and $r$,
$$
 \DD(C_n^r)=r(n-1)+1.
$$
\end{conjecture}
Although it was proved that $\mathsf D (G) > \mathsf D^* (G)$ for various classes of groups (e.g., \cite{Ge-Li-Ph12, Sa-Ch14a, Li20a}), the precise value of $\mathsf D (C_n^r)$ remained open (as late as 2018 the value of $\DD(C_6^4)$ was unknown
\cite[p.~606]{girard}).
In spite of such results,  Conjecture \ref{conj:homocyclic} has been addressed since decades (e.g., \cite{Ba-Sc80}). It developed into a significant,   long-standing problem in zero-sum
theory, and there was hope for an affirmative answer (see
\cite[Conjecture~A.5]{afk}, \cite[Conjecture~3.5]{gg06},
\cite[Conjecture~1]{girard}, and \cite[Conjecture~3]{wang2026}).
Also several recent results seemed to support this conjecture (e.g., \cite{Zh28a}).

In this paper, we disprove Conjecture~\ref{conj:homocyclic} for
all exponents $n>1$ that are not prime powers.

\Needspace{7\baselineskip}
\begingroup
\begin{theorem}\label{thm:large-rank}
Let $n>1$ have at least two distinct prime divisors. For any
factorization $n=PQ$ with $P,Q>1$ and $\gcd(P,Q)=1$, put
$$
 R_n=9(2n+1-P-Q).
$$
Then for every $r\ge1$,
$$
 \DD(C_n^r)\ge r(n-1)+1+
       \left\lfloor\frac{r-1}{R_n-1}\right\rfloor .
$$
In particular, $\DD(C_n^r)>r(n-1)+1$ for every $r\ge R_n$.
\end{theorem}
\par\endgroup
Note that Theorem \ref{thm:large-rank} also disproves a conjectured result on coverings (see \cite[Theorem 6.6]{Ga-Ge03a}).

Our second main result considers the behavior of $\mathsf D (C_n^r)$, when the rank is fixed and  the exponent tends to
infinity. This setting has been studied before. Indeed,  Girard
 proved that,  for every fixed
rank $r\ge1$,  $\DD(C_n^r)/n\to r$ as $n\to\infty$ (\cite[Theorem~1]{girard}).
In rank three, Pop, Rad and \c{T}urca\c{s} \cite[Corollary~C]{pop-rad-turcas}
used integral zero-sum incidence identities to prove that
$\DD(C_n^3)=3n-2$ whenever $n=mp^a$, where $m\in\{2,3,4\}$,
$p$ is a prime with $p\nmid m$, and $a$ is a nonnegative integer.

\begin{theorem}\label{thm:main}
Let  $r\ge3$ be a positive integer. Then there exists a constant $N(r)$ such that
for every  $n>N(r)$,
$$
                         \DD(C_n^r)=r(n-1)+1 .
$$
\end{theorem}
Our proof does not give an explicit value for $N(r)$.

\section{Notation and terminology}\label{sec:notation}
We follow notation and terminology used in \cite[\S2]{gg06}. Groups are written
additively. We denote by $C_n$ the cyclic group of order $n$, by $C_n^r$
the direct sum of $r$ copies of $C_n$, and by $\langle g\rangle$ the
subgroup generated by $g$.

We consider sequences over $G$ as  elements of the free abelian monoid
$\mathcal F(G)$ with basis $G$. We write $S=g_1\cdot\ldots\cdot g_\ell
=\prod_{g\in G}g^{\mathsf v_g(S)}$, where $\mathsf v_g(S)$ is the
multiplicity of $g$ in $S$. The element $1 = 1_{\mathcal F (G)}$ is called the {\it trivial sequence}. We write $R\mid S$
if $R$ is a subsequence of $S$, that is, if $\mathsf v_g(R)\le\mathsf v_g(S)$
for every $g\in G$. The length of $S$ is $|S|=\ell$ and its sum is
$\sigma(S)=g_1+\cdots+g_\ell$, with $\sigma(1)=0$. In particular,
$g^m$ denotes a sequence of $m$ copies of $g$, whose sum is $mg$.

The set of subsums is $\Sigma(S)=\{\,\sigma(R):R\mid S,\ R\ne1\,\}$.
We call $S$ a \emph{zero-sum sequence} if $\sigma(S)=0$, and
\emph{zero-sumfree} if $0\notin\Sigma(S)$. The maximal length of a
zero-sumfree sequence over $G$ is $\mathsf d(G)$, so
$\DD(G)=1+\mathsf d(G)$. For the invariant factors
$1<n_1\mid\cdots\mid n_r$ of $G$, put $\mathsf d^*(G)=\sum_{i=1}^r(n_i-1)$.
Taking $n_i-1$ copies of a generator of each cyclic summand gives a
zero-sumfree sequence, whence $\mathsf d(G)\ge\mathsf d^*(G)$. In particular, $\mathsf d^*(C_n^r)=r(n-1)$.

We identify $C_n^r$ with $(\mathbb Z/n\mathbb Z)^r$.
The standard basis element $e_i$ has a $1$ in position $i$
and zeros elsewhere. Thus $(x_1,\ldots,x_r)=\sum_{i=1}^r x_i e_i$, with
coefficients taken modulo $n$.
{For a prime $p$, $\mathbb F_p=\mathbb Z/p\mathbb Z$ is the field of $p$ elements.}
For integers $a\le b$, let $[a,b]=\{a,a+1,\ldots,b\}$.
For a statement $P$, the
indicator $\mathbf{1}_P$ is $1$ if $P$ holds and $0$ otherwise.

\begingroup
\setlength{\emergencystretch}{1em}
Write
$h(S)=\max_g\mathsf v_g(S)$ and
$\operatorname{supp}(S)=\{g:\mathsf v_g(S)>0\}$.
For $T\mid S$, the complementary subsequence is $T^{-1}S$.
A homomorphism acts termwise, preserving all indexed positions.
For $j\ge1$, put $\Sigma_j(S)=\{\sigma(T):T\mid S,\ |T|=j\}$.
Then $\Sigma(S)=\bigcup_{j\ge1}\Sigma_j(S)$. {We reserve $\Sigma_0(S)=\Sigma(S)\cup\{0\}$ for sums of subsequences of arbitrary length, including the empty subsequence.}

In a cyclic coordinate, $[a,b]$ denotes the image of the integer
interval modulo $n$ and is empty for $a>b$. {For a residue $x$ modulo $n$, $[x]_n$ denotes its representative in $[0,n-1]$.} For $x\in C_n$, a signed
representative attains ${\lVert x\rVert_n}=\min_{t\equiv x\bmod n}|t|$, with
either choice allowed at a tie.

For an ordered basis $(e_1,\ldots,e_r)$ of $C_n^r$, write $H_0=\{0\}$
and $H_j=\langle e_1,\ldots,e_j\rangle$ for $1\le j\le r$.
These subgroups form a flag. Its $j$th \emph{ordinary row} is the coset
$e_j+H_{j-1}$. Its population in $S$ counts indexed positions, including
all repetitions. {The \emph{stock} of a row consists of the indexed positions retained in it. The \emph{exterior} consists of all positions outside the chosen stocks, including any excess positions in a row.} Write $T_r(n)=r(n-1)+1$ and call a sequence of this
length \emph{critical}.

Asymptotic statements in Sections~\ref{sec:geometry}--\ref{sec:positions} concern
families with $n\to\infty$. In arguments by contradiction we may
pass to subsequences of exponents. {For $g>0$, write $f=O(g)$ when $|f|/g$ stays bounded and $f=o(g)$ when $f/g\to0$. For $f\ge0$, write $f=\Theta(g)$ when $f/g$ stays between two positive constants. A subscript on an asymptotic bound indicates a parameter on which its constants may depend.}
For nonnegative quantities, $f\gg g$ means $f\ge c g$ for some constant $c>0$ independent of $n$.
For $k>0$ a quantity is $\omega(k)$ if its ratio to $k$ tends to
infinity. In geometric arguments $[0,\lambda]v$ denotes the real line
segment $\{tv:0\le t\le\lambda\}$.
\par\endgroup

\section{Counterexamples in large rank}\label{sec:large-rank}

\begingroup
If $n>1$ has at least two distinct prime divisors, it admits
a factorization $n=PQ$ with $P,Q>1$ coprime. For example,
take $P$ to be the largest power of a chosen prime divisor
of $n$ that divides $n$. Fix any such factorization and put $B=n-1$.
For a vector $X$ in a power of $C_n$, let $X_P,X_Q$ denote
its two components under the Chinese remainder theorem. We regard
these as vectors in the same group, so that $X=X_P+X_Q$.
Write $e$ for the integer in $[0,n-1]$ congruent to $1$ modulo
$P$ and $0$ modulo $Q$. Thus $X_P=eX,\ X_Q=(1-e)X$.
We use a bit, equal to zero or one, to indicate whether a given
indexed term belongs to a subsequence.

\begin{lemma}\label{lem:two-choices}
Let $c,x$ be the bits of two distinct terms, and let $U,V$
be vectors with $PU=QV=0$. By adjoining $B$ terms and one
coordinate and adding coefficients to the two terms, one can
obtain the following additional contributions on the old coordinates
when the new coordinate sums to zero. If $cx=0$ the contribution
is zero, and if $cx=1$ it is either $U$ or $V$, with both
choices possible. At $c=x=0$ the auxiliary selection is empty.
\end{lemma}

\begin{proof}
Put $K=e-1$ and $D=U-V$. Since $1<e<n$, we have $1\le K<B$.
Give both terms coefficient one on the new coordinate and add $KD$
to each on the old coordinates. Use $B$ units on the new coordinate,
$K$ of them carrying $-D$. For $c+x=0,1,2$ one must take
respectively none, all, or all but one of these units. The first two
cases give zero. In the last case the contribution is $KD=V$
or $(K+1)D=U$.
\end{proof}

\begin{lemma}\label{lem:uniform-crt-switch}
Suppose $W$ is a zero-sumfree sequence of length $sB$ over
$C_n^s$. For $m\ge1$ specified terms $X_i$ and a permutation
$\pi$ of their indices, form $W'$ by replacing them by
$$
 A_i=(X_i)_P+(X_{\pi(i)})_Q .
$$
If a point $t_*$ is absent from $\Sigma_0(W')$, there is a
zero-sumfree sequence of length $RB+1$ in rank $R=s+2m$,
containing $n-2$ copies of a basis element.
\end{lemma}

\begin{proof}
Add one term of value $-t_*$ with bit $z$. Introduce $m$ new
coordinate units $f_i$, placing $B$ copies of each in the sequence.
Give one of the copies of $f_i$ the additional value $-A_i$.
Apply Lemma~\ref{lem:two-choices} to $z$ and the bit of $X_i$,
with alternatives
$$
 (X_{\pi^{-1}(i)}-X_i+f_{\pi^{-1}(i)})_P,\qquad
 (X_{\pi(i)}-X_i+f_i)_Q ,
$$
using a separate coordinate for each $i$.
If $z=0$, all added contributions vanish, the $f_i$-terms must
be absent, and zero-sumfreeness of $W$ forces the empty selection,
also on every auxiliary coordinate.

Suppose $z=1$. Each selected $X_i$, together with its correction,
gives either $A_{\pi^{-1}(i)}+e f_{\pi^{-1}(i)}$ or
$A_i+(1-e)f_i$. Each index $j$ occurs at most twice, once
with coefficient $e$ and once with coefficient $1-e$ on $f_j$.
If it occurs twice, these coefficients sum to one. All $B$ terms
on $f_j$ must then be used, cancelling exactly one copy of $A_j$.
If $j$ occurs once, at most one copy is cancelled. If it does
not occur, no term on $f_j$ is selected. Thus the old-coordinate
contribution apart from $-t_*$ belongs to $\Sigma_0(W')$,
a contradiction. We have added $2mB+1$ terms and $2m$
coordinates. Each $f_i$ has $n-2$ pure copies.
\end{proof}

\begin{proposition}\label{prop:separate-endpoints}
For each $m\in\{P,Q\}$, suppose $T^{(m)}$ is a zero-sumfree
sequence of length $h_m B$ over $C_n^{h_m}$. Distinguish
$b_m\ge1$ positions in $T^{(m)}$. {Let $\widetilde T^{(m)}$ be obtained from $T^{(m)}$ by projecting the distinguished terms onto the component of order $n/m$ and leaving all other terms unchanged. Assume $t^{(m)}\notin\Sigma_0(\widetilde T^{(m)})$.}
Put $h=h_P+h_Q,\ b=b_P+b_Q$ and
\begin{equation}\label{eq:disjoint-rank}
 R=h+9b-11 .
\end{equation}
Then $\DD(C_n^r)\ge rB+1+\lfloor(r-1)/(R-1)\rfloor$ for all $r\ge1$.
\end{proposition}

\begin{proof}
We construct $W'$ first. In each of two arrays a selection attaining
the prescribed target will have to select at least one initial term.
The operations in the array force any such selection to contribute
one at its last edge. The two arrays share this edge, where the target
excludes the simultaneous contributions. We then rearrange components
to obtain a zero-sumfree $W$ and apply
Lemma~\ref{lem:uniform-crt-switch}.

\emph{The two arrays.}
Use $h$ coordinates for the direct sum of the two groups supporting
$T^{(P)},T^{(Q)}$. Call these the original coordinates. Include
all undistinguished terms once. In the left array put $m=P$,
and in the right array put $m=Q$. In each array number the
distinguished terms $D_1,\ldots,D_{b_m}$. Start row $i$ with
a term called its endpoint, having selection bit $c_i$.
On the original coordinates this endpoint carries $(D_i)_m$.
Place the two complementary components $(D_1)_{n/m}$, one
from each array, together in a single additional term.

In each array, for $i=2,\ldots,b_m$ in order, place an operation
on rows $1,i$. It has two incoming edges, one from the preceding
vertex of row one and one from endpoint $i$. The first vertex of
row one is endpoint one. Join the last vertices of row one in the
two arrays by a central edge. When $b_m=1$, the last vertex is
the endpoint itself.

Each edge has a separate $C_n$ coordinate and $n-2$ pure copies
of its generator. On each incoming edge record the output of the
preceding vertex modulo $n/m$ and the input of the operation
modulo $m$. The central edge records one output from each array.
Thus on every edge the two specified residues determine a single
value in $C_n$. The output of an endpoint is its selection bit.

The operation on rows $1,i$ uses two terms with selection bits
$x_i,y_i$. The term with bit $x_i$ carries $(D_i)_{n/m}$ on
the original coordinates. Before any correction, the inputs
$S_1,S_i$ on the incoming edges and the output $U_1$ on the
next edge of row one are
$$
 S_1=y_i,\qquad S_i=x_i-y_i,\qquad U_1=y_i .
$$
Apply Lemma~\ref{lem:two-choices} to $c_i,x_i$ with alternatives
$-1$ at $S_i$ or $1$ at $U_1$, each zero elsewhere.
Order the alternatives so that the first is supported modulo $P$.
Use a fresh auxiliary coordinate for every operation. All
coefficients left unspecified are zero. This defines $W'$.

At zero on the auxiliary coordinate the corrected values satisfy
\begin{equation}\label{eq:inner-directions}
 \begin{aligned}
 c_i=0,\ S_1\ne0 &\ \Longrightarrow\ U_1=1,\\
 c_i=1,\ S_i\ne0 &\ \Longrightarrow\ U_1=1 .
 \end{aligned}
\end{equation}
The first implication follows from $S_1=U_1=y_i$.
For the second, the possible pairs $(S_i,U_1)$ are
$(-y_i,y_i)$ and $(1-y_i,1+y_i)$. In either pair a nonzero
first entry forces the second to be one.

\emph{An omitted subsum of $W'$.}
Prescribe $(t^{(P)},t^{(Q)})$ on the original coordinates and
zero on the auxiliary coordinates. Without an endpoint in array
$m$, its original coordinates could only realize a subsum of the
projected sequence. Thus at least one endpoint in each array is
required. On each incoming edge forbid output one followed by input
zero. On the central edge forbid two outputs equal to one. To forbid
a pair, take its value in $C_n$ minus one as the target on the edge.
Indeed the $n-2$ units can supply every residue except $-1$.

In either array, let $i$ be the largest index with $c_i=1$.
If $i\ge2$, the forbidden pair forces $S_i\ne0$, so $U_1=1$
at operation $i$. If $i=1$, row one already starts with output
one. Each subsequent operation has $c_j=0$, and the incoming
edge along row one forces $S_1\ne0$. By
\eqref{eq:inner-directions} the output remains one. Hence both
arrays contribute one at the central edge, contradicting its target.

\emph{The zero-sumfree sequence $W$.}
Split each endpoint, each of the two terms at each operation and the
single additional term into their two components, including the
coefficients added by Lemma~\ref{lem:two-choices}. In array $m$
call the components modulo $m$ and $n/m$ the input and output
components, respectively. Form new terms by adding components in
the following pairs.
\begin{enumerate}
\item\label{pair:same-term} Within each array pair the input of endpoint $i$ with the
component carrying $(D_i)_{n/m}$. This is part of the additional
term for $i=1$, and the output of $x_i$ otherwise.
\item\label{pair:endpoint} For $i\ge2$ pair the output of endpoint $i$ with the input
of $x_i$.
\item\label{pair:chain} Along row one pair consecutive terms of the list
$c_1,y_2,\ldots,y_{b_m}$, using the output of one and the input
of the next. Pair the two final outputs across the central edge.
\end{enumerate}
Every component occurs exactly once, and each pair uses both moduli.
Keeping the $P$-components in place, a permutation of the
$Q$-components therefore recovers $W'$ from $W$.

In a zero sum of $W$, the original coordinates give a zero sum
in each $T^{(m)}$. All undistinguished terms and all pairs
in {item~\ref{pair:same-term}} are therefore absent. For each operation, the remaining
components of $c_i,x_i$ form the single pair in {item~\ref{pair:endpoint}}.
Consequently the two residues of $c_i+x_i$ in the calculation
of Lemma~\ref{lem:two-choices} are now both the bit of this pair.
Both the new coordinate and the added old coefficients depend
only on this sum. Since the bit is zero or one, the net correction
is zero. Each edge of row one now receives the two components of
a single bit by {item~\ref{pair:chain}}. The $n-2$ pure units cannot cancel one,
so all these bits vanish. In particular all $y_i$ are absent.
On each remaining edge the same argument excludes the pair in {item~\ref{pair:endpoint}}.
All pure units and auxiliary selections are then empty.

\emph{Rank and length.}
There are $d=b-2$ operations, $1+2d$ edges and $M=b+2d+1$
paired positions. The rank is $s=h+1+3d$. The length is
$$
 hB-b+M+(1+2d)(B-1)+dB=sB .
$$
Lemma~\ref{lem:uniform-crt-switch} gives a zero-sumfree sequence
of length $RB+1$ in rank $R=s+2M=h+9b-11$.

Delete its $n-2$ specified copies of a basis element $f$.
Every nonempty subsum of the remainder in $\langle f\rangle$
equals $f$, since the removed terms can cancel every other
multiple of $f$. Replace any term $a$ of a zero-sumfree
sequence by this remainder, sending $f$ to $a$ and using
$R-1$ fresh coordinates for the other basis elements. A zero
using any inserted terms would contract to a zero using $a$.
This replacement adds $R-1$ coordinates and $(R-1)B+1$
terms. Start with $B$ copies of a cyclic generator, repeat
$\lfloor(r-1)/(R-1)\rfloor$ times, and fill each spare
coordinate with $B$ copies of its generator.
\end{proof}

\bigskip

\Needspace{6\baselineskip}Now we are in a position to prove Theorem~\ref{thm:large-rank}.

\begin{proof}[Proof of Theorem~\ref{thm:large-rank}]
For each $m\in\{P,Q\}$ take $B$ copies of a generator
of $C_n$, distinguishing $b_m=n-m+1$ of them. The sequence
is zero-sumfree since $B<n$. After projecting the distinguished
terms onto the component of order $n/m$, only the $m-2$
undistinguished terms contribute modulo $m$. Thus the negative
of the generator is an omitted subsum.
Proposition~\ref{prop:separate-endpoints} applies with
$h=2$ and $b=2n-P-Q+2$. It gives the stated bound with
$R=h+9b-11=9(2n+1-P-Q)=R_n$.
For $r\ge R_n$ the floor is at least one, proving strict inequality.
\end{proof}

\begin{remark}
For particular exponents, the lower bound can be improved.
If $n=2k$ with $k\ge3$ odd, then $\DD(C_{2k}^r)\ge r(2k-1)+1+
       \left\lfloor\frac{r-1}{5k-1}\right\rfloor$ for every $r\ge 1$.
See Theorem~1.1 of \href{https://arxiv.org/abs/2611.24775v3}{arXiv:2611.24775v3}.
For $n=6$, we have
$\DD(C_6^r)\ge 5r+1+\lfloor 2(r-1)/15\rfloor$ for every $r\ge1$; We omit the construction here to keep
the exposition focused on the main theorems.
\end{remark}
\par\endgroup

\begingroup
\setlength{\emergencystretch}{2em}
\newcommand{\Z}{\mathbb Z}\newcommand{\R}{\mathbb R}\newcommand{\Q}{\mathbb Q}
\newcommand{\F}{\mathbb F}\newcommand{\rk}{\operatorname{rank}}
\newcommand{\bin}[2]{\binom{#1}{#2}}\newcommand{\T}{\mathbb T}
\newcommand{\norm}[1]{\lVert #1\rVert}

\section{Lattice selection and finite projected support}\label{sec:geometry}
In this section we prove two geometric results used in the proof of
Theorem~\ref{thm:main}. Theorem~\ref{geom-main} gives a strict selection
on a nonempty subset of a balanced list of lattice vectors.
Theorem~\ref{chart-main} produces a homomorphism under which all but
$o(n)$ terms of a near-extremal zero-sumfree sequence take values in
a fixed finite set. Section~\ref{sec:stability} uses these results to
construct ordinary flags.

\subsection{Lattice selection and covering}
Let $v_1,\ldots,v_s$ be nonzero vectors generating a lattice
$\Gamma$ of rank $r\ge2$ in its real span. Assign them positive real
capacities $\lambda_i$ with $\sum_i\lambda_i=r$. Suppose that in
every nonzero proper subspace spanned by columns, the total capacity
of the columns it contains is strictly less than its dimension.

\begin{theorem}\label{geom-main}
There exist a nonempty $J\subseteq\{1,\ldots,s\}$, real
coefficients $0<c_i<\lambda_i$ for $i\in J$, and integers $n_i$
for $i\in J$ such that
$$
                    \sum_{i\in J} c_i v_i=\sum_{i\in J} n_i v_i .
$$
\end{theorem}

\begin{lemma}\label{geom-metric}
In any fixed linear coordinates on the span there are
$\rho_i>0$ and a positive definite matrix $Q$ satisfying
$$
       \sum_i\rho_i v_i v_i^\top=Q^{-1},
       \qquad \rho_i v_i^\top Qv_i=\lambda_i .
$$
\end{lemma}
\begin{proof}
Let $V$ be the matrix with columns $v_i$ and $\mathcal B$ the
finite collection of index sets of column bases. {Write $V_I$ for the submatrix formed by the columns indexed by $I$.} For any nonempty
proper index set $I$, one has $\lambda(I)=\sum_{i\in I}\lambda_i
<\rk V_I$. For a nonspanning $I$ this follows by applying strict
balance on its span. For a spanning proper $I$ it follows from
the positive capacity of an omitted column and total mass $r$.

For any vector $x\in\R^s$ whose entries are not all equal, order
them as $x_{i_1}\ge\cdots\ge x_{i_s}$. Let
$I_j=\{i_1,\ldots,i_j\}$. Scan this order, retaining a column
exactly when it increases the span of the retained columns.
The resulting basis $B$ satisfies $|B\cap I_j|=\rk V_{I_j}$,
since after each scan its retained columns span that prefix.
Subtracting the two weighted sums and telescoping gives
$$
 \sum_{i\in B}x_i-\sum_i\lambda_i x_i
   =\sum_{j=1}^{s-1}(x_{i_j}-x_{i_{j+1}})
                        \big(\rk V_{I_j}-\lambda(I_j)\big)>0 .
$$
The contribution of $x_{i_s}$ cancels because both masses are $r$,
and at least one consecutive gap is positive. The maximum of the
left-hand expression over $B\in\mathcal B$ is continuous and
positive on the unit sphere of $\sum_i x_i=0$. By compactness
it has a positive lower bound $\eta$ there.

Consider the function
$$
 f(x)=\log\det\!\left(\sum_i e^{x_i}v_i v_i^\top\right)
                       -\sum_i\lambda_i x_i .
$$
The matrix is positive definite because the columns span. By
Cauchy--Binet its determinant is
$\sum_{B\in\mathcal B}(\det V_B)^2\exp(\sum_{i\in B}x_i)$.
Let $c$ be the minimum of ${\log\bigl((\det V_B)^2\bigr)}$ over the bases.
Thus $f(x)\ge\eta\|x\|+c$ on $\sum x_i=0$ and $f$ attains a
minimum on this hyperplane. It is unchanged
by adding the same constant to every coordinate, so this is also
an unconstrained minimum. Its partial derivatives vanish. Set
$\rho_i=e^{x_i}$ at the minimum and
$Q=(\sum_i\rho_i v_i v_i^\top)^{-1}$. Differentiating the
log determinant gives
$0=\operatorname{tr}(Q\rho_i v_i v_i^\top)-\lambda_i
=\rho_i v_i^\top Qv_i-\lambda_i$ for every $i$.
\end{proof}

\begin{proof}[Proof of Theorem~\ref{geom-main}]
{Take $Q$ and the weights $\rho_i$ from Lemma~\ref{geom-metric}.} Use the inner product $\langle y,w\rangle=y^\top Qw$ and its norm.
Set $h=\frac12\sum_i\lambda_i v_i$ and let $P$ be the nonempty
finite set of sites in $\Gamma$ nearest to $h$. For any $x\in P$
put
$$
              t_i(x)=\lambda_i/2+\rho_i\langle x-h,v_i\rangle .
$$
The first identity of Lemma~\ref{geom-metric} gives
$\sum_i t_i(x)v_i=h+Q^{-1}Q(x-h)=x$. Comparing the squared
distance at $x$ with those at $x\pm v_i$ gives
$|\langle x-h,v_i\rangle|\le\|v_i\|^2/2$. The second identity
of the lemma now gives $0\le t_i(x)\le\lambda_i$. Equality in
the distance comparison yields
$$
    t_i(x)=0\ \Longrightarrow\ x+v_i\in P,\qquad
    t_i(x)=\lambda_i\ \Longrightarrow\ x-v_i\in P .
$$

Choose $x\in P$ and an integer representation $x=\sum_i n_i v_i$
minimizing $\sum_i|n_i|$ over all such pairs. Every site has a
representation since the columns generate $\Gamma$, so this
nonnegative integer has a minimum. If $t_i(x)=0$ and $n_i<0$,
then $x+v_i\in P$ has the representation obtained by replacing
$n_i$ by $n_i+1$. Its cost is smaller by one, a contradiction.
Similarly if $t_i(x)=\lambda_i$ and $n_i>0$, use $x-v_i$ and
$n_i-1$. Thus
$$
          t_i(x)=0\ \Longrightarrow\ n_i\ge0,\qquad
          t_i(x)=\lambda_i\ \Longrightarrow\ n_i\le0 .
$$

Take $c_i=(1-\epsilon)t_i(x)+\epsilon n_i$ with
$0<\epsilon<1$ sufficiently small. Both representations have
sum $x$, hence so does this one. If $t_i(x)$ is strictly interior,
then $c_i$ remains strictly interior by continuity. At an upper
endpoint $n_i\le0<\lambda_i$, so $c_i<\lambda_i$ and it remains
positive for sufficiently small $\epsilon$. At a lower endpoint,
either $n_i=0$ and $c_i=0$, or $n_i>0$ and
$0<c_i=\epsilon n_i<\lambda_i$ for small $\epsilon$. A single
choice of $\epsilon$ works for the finite list. In particular
$c_i=0$ only when $n_i=0$. On $J=\{i:c_i>0\}$ we therefore have
$x=\sum_{i\in J}c_iv_i=\sum_{i\in J}n_iv_i$.

It remains to check that $J$ is nonempty. Otherwise all
$t_i(x)=n_i=0$, whence $x=0$, $0\in P$, and every $v_i\in P$.
In particular $2\langle h,v_i\rangle=\|v_i\|^2>0$ for every $i$.
The generated cone is pointed because $h$ is strictly positive
on its generators. Choose a column $u$ on an extreme ray. Such
a ray exists: intersect the cone with $\langle h,v\rangle=1$ and
choose a vertex of the convex hull of the normalized columns.
Every column on this ray equals $u$, since the same displayed
equality for $tu$ and $u$ implies $t=t^2$ for $t>0$.
In the nonnegative reconstruction of the nearest site $u$,
all positive terms lie on this extreme ray. For a column equal
to $u$, its coefficient at $u$ is
$t_i(0)+\rho_i\|u\|^2=\lambda_i$. Therefore
$u=(\sum_{i:v_i=u}\lambda_i)u$, and its ray has total capacity
one. No column lies on the opposite ray because $\langle h,v_i\rangle>0$.
The line through $u$ is proper for $r\ge2$, contradicting strict
balance. This proves $J\ne\varnothing$.
\end{proof}

\medskip\noindent\textit{Covering every translate with strict surplus.}\quad
For an
integer matrix $B$ with $b$ independent rows {and $s$ columns}, let $w_i$ be its columns
and $B_I$ the submatrix on $I$. Set $\T=\R/\Z$.
\begin{lemma}\label{cov-cover}
If $\lambda_i>0$ and $\sum_{i\notin I}\lambda_i>b-\rk B_I$
whenever $\rk B_I<b$, then each coset of
$\ker(B:\T^s\to\T^b)$ meets the image of
$\prod_i(0,\lambda_i)$. The distances of the chosen real coordinates
from their respective endpoints can be bounded below by a positive
constant independent of the coset.
\end{lemma}
\begin{proof}
The case $b=0$ follows by taking each coordinate at the midpoint.
Suppose $b>0$. Choose $\delta>0$ so that $\mu_i=\lambda_i-2\delta>0$
and $\sum_{i\notin I}\mu_i\ge b-\rk B_I$ for all $I$.
This is possible by strictness and finiteness; the inequality for a
spanning $I$ follows from nonnegativity.
Let $P$ be the convex hull of the indicator vectors $u_J$ of the
column bases. There is a $\tau'\in P$ with $\tau'_i\le\mu_i$:
otherwise a separating linear functional for $\mu$ and
$P+\R_{\ge0}^s$ gives a cost vector $c\ge0$ with
$c\cdot\mu<\min_J c\cdot u_J$. To see that this is impossible, order
the columns by increasing cost, retaining a column when it increases
the rank of those retained. Below each cost threshold $t$ the retained
columns span all columns of cost $<t$. Thus the selected basis has
$b-\rk B_{\{i:c_i<t\}}$ terms of cost at least $t$, up to the finitely
many threshold values. Integrate in $t$ from zero to $\max c_i$.
The capacity inequalities give $c\cdot\mu$ at least the cost of this basis.

Mix $\tau'$ with a sufficiently small positive proportion of the average
of all the $u_J$. The resulting $\tau$ is in the relative interior of
$P$ and has $\tau_i<\lambda_i-\delta$. There are positive numbers
$\rho_i$ such that, with $Q=(\sum_i\rho_iw_iw_i^{\mathsf T})^{-1}$,
$$
                      \rho_i w_i^{\mathsf T}Qw_i=\tau_i .
$$
Indeed, set $\rho_i=e^{x_i}$ and expand the determinant in
$$
 f(x)=\log\det\Big(\sum_i e^{x_i}w_iw_i^{\mathsf T}\Big)-\tau\cdot x
     =\log\sum_J(\det B_J)^2e^{(u_J-\tau)\cdot x}
$$
by Cauchy--Binet. On the span of the $u_J-\tau$, the maximum of
$(u_J-\tau)\cdot x$ on any unit direction is positive because $\tau$
is in the relative interior. It has a positive uniform lower bound
by compactness, unless that span is zero. Hence $f$ attains a minimum
there, also in the zero-dimensional case. It is constant along the
orthogonal complement. All its partial derivatives vanish at the
minimum, and differentiating the log determinant gives the assertion.

For any $h\in\R^b$ take $z\in\Z^b$ minimizing
$(z-h)^{\mathsf T}Q(z-h)$. Comparing with $z+w_i$ and $z-w_i$ gives
$|w_i^{\mathsf T}Q(z-h)|\le w_i^{\mathsf T}Qw_i/2$. Thus
$$
 z-h=\sum_i \rho_i w_i w_i^{\mathsf T}Q(z-h)
     \in B\prod_i[-\tau_i/2,\tau_i/2].
$$
Given a desired real output $y$ modulo $\Z^b$, put
$h=B(\lambda/2)-y$. The last display gives
$Bx=y+z$ for some $x_i\in[(\lambda_i-\tau_i)/2,(\lambda_i+\tau_i)/2]$.
These intervals lie at distance at least $\delta/2$ from the endpoints
of the original box. This proves the claim.
\end{proof}

\subsection{Finite projected support at arbitrary exponents}\label{sec:projection}

Let $n\to\infty$, fix $d\ge1$, and put $G_n=C_n^d\oplus E_n$,
where $|E_n|\le C$ {for a fixed constant $C$} and $nE_n=0$. A sequence is zero-sumfree if none
of its nonempty subsequences has sum zero. All assertions about limits
below allow passage to a subsequence. The terms may all be different.

\begin{theorem}\label{chart-main}
Let $S_n$ be zero-sumfree in $G_n$ with $|S_n|/n\to d$. {After passing to a subsequence, there} exist
$1\le k\le d$ and homomorphisms $F_n:G_n\to C_n^k$ of uniformly
bounded image index with the following properties.
\begin{enumerate}
\item The kernel is isomorphic to $C_n^{d-k}\oplus E'_n$ with
$|E'_n|$ bounded. It contains $(d-k)n+o(n)$ terms.
\item Except for $o(n)$ positions the outputs are reductions modulo
$n$ of one fixed finite set in $\Z^k$. Its nonzero members decompose
into families in complementary rational subspaces $W_1,\ldots,W_a$.
The limiting number of terms in family $j$, divided by $n$, equals
$t_j=\dim W_j$. In every nonzero proper subspace of $W_j$, the
limiting mass of that family is strictly less than the dimension.
\item In a positive definite metric on $W_j$, each of its columns is
a difference of two lattice points nearest a common point in
$W_j$ for the lattice $W_j\cap\Z^k$. Consequently the family has at
most $2^{t_j}(2^{t_j}-1)$ distinct nonzero columns.
\end{enumerate}
When $k=1$, after a choice of sign, $n+o(n)$ positions have output $1$
and only $o(n)$ positions have output outside $\{0,1\}$. {In this one-dimensional case},
$F_n$ is surjective for all sufficiently large $n$ {on this subsequence}.
\end{theorem}

{The proof of Theorem~\ref{chart-main} follows Lemma~\ref{chart-f:compact}.}

\medskip\noindent\textit{{A subgroup index bound.}}\quad
For a finite {abelian} group $H$ {denote by $\exp H$ the least positive integer $m$ such that $mH=0$, and} let $\eta(H)$ be the least length guaranteeing
a nonempty zero sum of length at most $\exp H$. We first establish
the subgroup estimate used to produce the projection.
\begin{lemma}\label{chart-macro}
Suppose $H_n$ has at most $R$ invariant factors, all dividing $n$.
Pass to a subsequence where their ratios to $n$ have limits, padding
by ones. If $j$ limits are positive, then
$\limsup \D(H_n)/n\le j$. In particular a zero-sumfree sequence
in $G_n$ of length $>(d-1+\gamma)n$ generates a subgroup of index
bounded in terms of $d,C,\gamma>0$.
\end{lemma}
\begin{proof}
Let $t$ be the largest invariant factor of ratio tending to zero,
or one if there is none. All smaller invariant factors divide $t$.
Embed $H_n$ in $C_n^j\oplus C_t^{R-j}$ by embedding each cyclic
factor in its corresponding coordinate. {Let $c_R$ be the constant depending only on $R$ in the Alon--Dubiner bound.} If $j=0$,
$\D(H_n)\le\eta(C_t^R)\le c_R t=o(n)$ by Alon and Dubiner~\cite{alon-dubiner}.
{Their result guarantees a zero sum of exactly $t$ terms, which suffices for the bound on $\eta$.}
For $j>0$ reduce the first $j$ coordinates modulo $t$. The quotient
is $C_t^R$ and its kernel is $C_{n/t}^j$. From a sequence of length
$(\D(C_{n/t}^j)-1)t+\eta(C_t^R)$ extract disjoint nonempty
quotient-zero blocks of length at most $t$. After the first
$\D(C_{n/t}^j)-1$ extractions at least $\eta(C_t^R)$ positions
remain, so one further block exists. Some nonempty selection of
their $\D(C_{n/t}^j)$ sums vanishes in the kernel. Since $t=o(n)$,
Girard's fixed-rank limit~\cite{girard} (or the cyclic formula for $j=1$)
gives $\D(C_{n/t}^j)t/n\to j$, proving the bound.

A subgroup of $G_n$ has bounded rank and order at most $C n^d$.
If its index tends to infinity then its order is $o(n^d)$,
since $|G_n|\le Cn^d$. It can have at most $d-1$ invariant factors
with positive limiting ratio, since $d$ such factors would already
give order at least a fixed positive multiple of $n^d$.
The first assertion contradicts the proposed zero-sumfree length.
This subsequence argument proves uniform boundedness of the index.
\end{proof}

\medskip\noindent\textit{Isolated exceptions to open lattice covering.}\quad
For a lattice $\Lambda$ in a space $V$ and segment sum $B$ spanning
$V$, call the missing targets the elements of
$V/\Lambda$ outside the image of $\operatorname{int}_V B$.
\begin{lemma}\label{chart-isolated}
Let the nonzero columns $v$ in a lattice $\Lambda$ span its real
space of dimension $t\ge1$. Give them positive masses $\lambda_v$
of total $t$ and finite first moment, with mass strictly less than
the dimension in every nonzero proper subspace. The missing targets
for $B=\sum[0,\lambda_v]v$ form a finite set. If any exist, {there are a positive definite metric and a point $x\in V$ such that each column is the difference of two points of $\Lambda$ nearest to $x$ in this metric.}
\end{lemma}
\begin{proof}
Put $c=\frac12\sum\lambda_vv$. There is a positive definite matrix
$Q$ satisfying, in orthonormal coordinates,
$$
 Q^{-1}=\sum_v\lambda_v vv^{\mathsf T}/(v^{\mathsf T}Qv).
$$
For $t=1$ use $Q=1$. In higher dimensions apply the
determinant-one scatter functional of D\"umbgen, Pauly and Schweizer,
{\cite[Theorems~4.9 and~6.1]{scatter}}, to the probability law placing mass
$\lambda_v/t$ at $vv^{\mathsf T}/\|v\|^2$. Each matrix is nonzero
of trace one, and their range probabilities on every nonzero proper
subspace are strictly less than its dimension divided by $t$.
These are the strict subspace hypotheses of their scale-invariant
case. {Writing $\Sigma$ for their determinant-one scatter matrix and $Q=\Sigma^{-1}$, the factor $t$ in their estimating equation cancels the factor $1/t$ in the masses. The trace normalization of each rank-one matrix cancels in the ratio $vv^{\mathsf T}/(v^{\mathsf T}Qv)$, giving the displayed identity.} The resulting inverse scatter matrix $Q$ is continuous under
small perturbations of these masses in total variation. Strictness
persists uniformly under sufficiently small perturbations: on the
compact Grassmannian of each dimension the subspace mass is upper
semicontinuous, by finite truncation of the measure.

If $y\bmod\Lambda$ is missing, take a nearest lattice point $z$
to $x=c-y$. Comparing $z$ with $z\pm v$ gives
$|v^{\mathsf T}Q(z-x)|\le v^{\mathsf T}Qv/2$. The displayed
identity expresses $z-x$ as a sum of the columns with coefficients
bounded by $\lambda_v/2$ in absolute value. Adding $c$ represents
$z+y$ in the closed $B$. For each fixed $v$ at least one nearest
point must give equality in this coordinate. Indeed otherwise
replace $\lambda$ by $(1-\epsilon)\lambda+\epsilon t\,\mathbf1_v${, where $\mathbf1_v$ denotes unit mass at $v$,} while keeping $x$ fixed. The metrics converge to $Q$. All new
nearest points eventually belong to the finite set of old nearest
points: the metrics are uniformly comparable, and on any bounded
lattice set the old distances above the minimum have a positive gap.
Thus at coordinate $v$ the reconstructed coefficient still lies
strictly inside the old segment, by the assumed strict inequality
and continuity. At every other coordinate the new half-mass is
strictly smaller than the old. After adding $c$, varying the
coefficients at $t$ independent columns gives a neighborhood of
$z+y$ in $B$. This contradicts the definition of a missing target.

Equality means $z$ and one of $z\pm v$ are both nearest to $x$.
The nearest points therefore span affinely, since their differences
include the spanning columns. Any $t+1$ affinely independent nearest
points determine $x$ uniquely by their squared-distance equalities:
subtracting the first equation from the others gives $t$ independent
linear equations in $Qx$. Move $x$ to a compact fundamental
parallelepiped modulo $\Lambda$. Its nearest points belong to one
fixed finite lattice set (compare with zero). Thus there are only
finitely many such $x$, proving the first assertion. In fact the
number of missing targets is at most $t!$. For each such $x$,
the convex hull of its nearest lattice points contains a lattice
simplex of volume at least $\operatorname{covol}(\Lambda)/t!$.
These full-dimensional hulls have disjoint interiors for distinct
centers, including lattice translates. To verify this, use the
{inner product $\langle x,u\rangle_Q=x^{\mathsf T}Qu$ and the} affine function $l_x(u)=2\langle x,u\rangle_Q+
\operatorname{dist}_Q(x,\Lambda)^2-\|x\|_Q^2$.
At lattice sites $u$, $l_x(u)\le\|u\|_Q^2$ with equality exactly
at the nearest sites. For two centers $x,x'$ the function
$l_x-l_{x'}$ is consequently nonnegative on the first hull and
nonpositive on the second. Overlap of their interiors makes it
zero on an open set and equality of the slopes gives $x=x'$.
Integrating volumes modulo $\Lambda$ now yields the asserted
uniform count. Distinct
nearest points have distinct residues modulo $2\Lambda$, since
otherwise their integral midpoint is strictly nearer. This also
gives {at most $2^t(2^t-1)$ oriented differences between distinct nearest points}.
\end{proof}

\subsection{Construction of the initial projection}\label{chart-exchange}
In this subsection write $r=d$, $G=C_n^r\oplus E_n$, and $I=I(r,C,1/2)$ from Lemma~\ref{chart-macro}. Write $I(1/2)$ for this bound. We have $n^r\le|G|\le Cn^r$ and $nG=0$. At a stage with $k$ independent output coordinates the kernel has order $O(n^{r-k})$, by bounded image index. These are the group-size inputs to the exchange argument.

\begin{lemma}\label{chart-l:column-limits}
Let $L_j$ be finite sequences in $\mathbb Z^q$, and let $h_j\to\infty$
with $|L_j|/h_j\to m$. Write $\mu_j(v)={\mathsf v_v}(L_j)/h_j$ and
assume $\sum_v\mu_j(v)\|v\|\le C$ for all $j$. On a subsequence there
are weights $\lambda_v\ge0$ such that for every $A\subseteq\mathbb Z^q$,
$$
 \sum_{v\in A}\mu_j(v)\longrightarrow\sum_{v\in A}\lambda_v,\qquad
 \sum_v\lambda_v=m,\qquad \sum_v\lambda_v\|v\|\le C.
$$
Let $Z$ be the closure of the union of the finite sums
$\sum_{v\in F}[0,\lambda_v]v$ over finite $F\subseteq\mathbb Z^q$.
Every point of $Z$ is a limit of subsums of $L_j$ divided by $h_j$, {with the empty subsum allowed}.
\end{lemma}
\begin{proof}
Since $\mu_j(v)\le |L_j|/h_j$, diagonal extraction on the countable
lattice gives limits $\lambda_v$ for all $v$. The moment bound passes
to each finite ball and then, by increasing the radius, to the whole
lattice. Thus both masses outside radius $R>0$ are at most $C/R$.
The sums on $A$ converge inside the ball, so the upper limit of
their absolute difference on $A$ is at most $2C/R$. Let $R\to\infty$. Taking
$A=\mathbb Z^q$ gives total mass $m$.

For $x=\sum_{v\in F}t_v v$ in a finite sum, take
$\min\{\lfloor h_jt_v\rfloor,{\mathsf v_v}(L_j)\}$ copies of each $v\in F$.
Their counts divided by $h_j$ tend to $t_v$ because $0\le t_v\le\lambda_v$.
Thus the normalized subsum tends to $x$. For $x$ in the closure,
first choose a finite-sum point within $\eta$ of $x$. The nearest
normalized subsum to $x$ then has distance at most $2\eta$ for all
large $j$. Let $\eta\downarrow0$.
\end{proof}

A useful consequence {of Lemma~\ref{chart-macro}} concerns any character $\chi{:G\to C_n}$ vanishing on more
than $(r-1+1/2)n$ terms of a zero-sumfree sequence. They generate a
subgroup inside its kernel, so $|\chi(G)|\le I(1/2)$. Conversely,
any nonzero representative $m$ of a value of $\chi$ gives
$|\chi(G)|\ge n/\gcd(n,m)\ge n/|m|$. We use these two inequalities
to check that the coordinate columns span at each stage.

{A convex body is hollow in a lattice if its relative interior contains no lattice point.}
For a {full-rank} lattice $\Gamma\subset\mathbb R^k$, write
$\Gamma^*=\{u:u\cdot\Gamma\subset\mathbb Z\}$. A class of positions
will be assigned a single column $a\in\mathbb Z^k$. To realize a
covector $u$ on these classes means that some character of $G$
takes the value $u\cdot a\bmod n$ at every position assigned $a$.

A character {$\chi$} of a subgroup $J$ {of $G$} extends to $G$. In circle notation,
choose $g\notin J$, let $t$ be the order of $g+J$, and take
$\xi^t=\chi(tg)$. Then $\chi(j+mg)=\chi(j)\xi^m$ is well-defined:
any two representations differ by an integer multiple of the relation
$tg\in J$. Repeating over generators extends $\chi$ to $G$, where
all values are $n$th roots of unity because $nG=0$.
For rounding counts, let $A$ be a fixed integer matrix, $L=\ker_{\mathbb Z}A$,
and $D\ge1$ a fixed integer. For any $x\in\mathbb R^s$ and
$b\in\mathbb Z^s$ with $Ax=Ab$, there exists
\begin{equation}\label{chart-eq:count-rounding}
         u\in b+DL,\qquad \|u-x\|_\infty\le C(A,D).
\end{equation}
Indeed $x-b\in\ker_{\mathbb R}A$. Row reduction over $\mathbb Q$
followed by clearing denominators gives a real basis $v_1,\ldots,v_e$
of that kernel in $L$. Write $x-b=\sum_i t_iDv_i$, and round each
$t_i$ to a nearest integer $m_i$. Then $u=b+\sum_i m_iDv_i$ belongs
to the required coset and
$\|u-x\|_\infty\le(D/2)\sum_i\|v_i\|_\infty$. The argument includes
the zero-dimensional kernel. In particular, suppose the count map
$g:\mathbb Z^s\to G$ gives the actual group sum, $g(b)=0$, and
$Dg(L)=0$. Then $g(u)=0$, since $u-b\in DL$. {Let $N_i$ be the available count in class $i$.
If $s\ge1$ and $C(A,D)<x_i<N_i-C(A,D)$ for every $i$,} $u$
gives a nonempty zero-sum subsequence. For fixed tolerance only
finitely many matrices occur, so the rounding error is uniform.

{An exchange replaces a selected group of positions (its first arm) by a disjoint group (its second arm). Its difference is the sum of the second arm minus the sum of the first. A pool is a collection of exchanges on disjoint positions. Its optional sums are the sums of differences over subsets of the pool, including the empty subset.}

\begin{lemma}\label{chart-f:chart}
Fix a sufficiently small $\epsilon>0$. {Let $(S_n)$ be a family of zero-sumfree sequences in {$G$} with $|S_n|/n\to r$.} For all sufficiently large $n$,
{$S_n$} has disjoint classes of sizes
$N_j$ that cover all but $\epsilon n$ positions, together with columns
$a_j$ generating $\mathbb Z^k$ for some $1\le k\le r$ and a lattice
$\Gamma\supseteq\mathbb Z^k$. Every $u\in\Gamma^*$ is realized on the
classes. The contraction about the center by factor $1-2\epsilon$ of
$Z=\sum_j[0,N_j/n]a_j$ is hollow in $\Gamma$. At a fixed tolerance
there are only finitely many choices of $k$, the list of columns, and
$\Gamma$, and $N_j/n$ is bounded below by a positive constant depending
only on the tolerance.
\end{lemma}
\begin{proof}
\smallskip\noindent\emph{The induction data.}
Start at $k=0$ with the trivial map and one class containing the whole
sequence. At stage $k$ a map $F:G\to C_n^k$ of bounded image index
takes a fixed bounded integer value $y_j$ on class $j$. All bounds
may depend on the tolerance $\epsilon$. There are at most $r+1$
stages, including zero. Allocate $\delta n$ losses per stage, where
$\delta=\epsilon/(r+1)$. Discard classes of size below $\tau n$.
Since the number of classes is bounded, choose $\tau>0$ small enough
that this costs at most $\delta n/4$. The total budget
$\epsilon n$ leaves more than $(r-1+1/2)n$ positions for large $n$,
since {$|S_n|/n\to r$} and $\epsilon<1/4$.
For $k>0$ the retained columns span $\mathbb Q^k$. Otherwise
a nonzero integral row $\ell$ annihilates them and can be chosen
from a finite list at this tolerance. Its image on $C_n^k$ has order
at least $n/\|\ell\|_\infty$, by evaluating on a basis vector.
Restriction to $\operatorname{im}F$ costs at most its bounded index.
Thus $|(\ell F)(G)|$ grows, contrary to the preceding character bound.
Choose a basis matrix $B$ of the integer lattice
generated by the $y_j$ and set $a_j=B^{-1}y_j$. The matrix
$A=(a_j)_j$ has integral columns generating $\mathbb Z^k$.
We either add a coordinate or find the contracted hollow body. The
following budget is fixed before testing any subgroup at this stage.
Choose $I=I(1/2)$ from Lemma~\ref{chart-macro}. Choose $C\ge1$ with
$|\ker F|\le C n^{r-k}$ for this chart, using its known image index.
Put
$$
 L=1+\lceil\log_2(CI)\rceil,\qquad
 \eta=\min\{\delta/(8L),\,\epsilon\tau/(32L)\}.
$$
At each successful count test give each of the two pools at most
$\eta n$ endpoints. Keep these endpoints out of the sample for the
character and discard them before a descent or augmentation. Give the
centering deletions together at most $\eta n$ per test. At each new
subgroup test the current count box afresh. During the first $L$
tests the total endpoints and deletions are at most $3L\eta n$,
less than $\delta n/2$ and $\epsilon\tau n/8$. Every initial
class therefore retains at least $\tau n/2$ positions,
and every retested coordinate has slack at least $\epsilon\tau n/2$
from zero and its current full count. The same budget bounds losses while the pools are under construction.
Below we use the retained length and the slack in each class to
show that the first $L$ tests suffice.
\smallskip\noindent\emph{The exchanges and the count test.}
Fix an integral right inverse $U$ of $A$ by writing each standard
basis vector as an integer combination of its columns.
An exchange between two terms of one class uses one position in
each arm. Also, for each nonzero $w_j=e_j-Ua_j$, use two disjoint
arms with respective class counts $w_j^+$ and $w_j^-$, the positive
and negative parts. Since $Aw_j=a_j-AUa_j=0$, their group sums have the same $F$-value. These exchanges use a bounded number of
positions because $A,U$ are fixed at this stage. At $k=0$, $w_j=e_j$,
so a singleton may be exchanged for the empty selection.
Let $K$ be generated by the differences of all such exchanges on
the unused positions. Each unused class is constant modulo $K$
by the two-term exchanges. Its label therefore defines a homomorphism
${\Phi}:\mathbb Z^{\#\text{classes}}\to G/K$. A representative exchange
for each $w_j$ exists because the unused classes still have linear
size. Hence ${\Phi}w_j=0$ for each $j$. For $t\in\ker_{\mathbb Z}A$,
$t=({\mathrm{Id}}-UA)t=\sum_jt_j w_j$, so ${\Phi}t=0$. Thus ${\Phi}$ factors through
$A$, and its image is generated by $k$ elements of exponent
dividing $n$. Its order is at most $n^k$. Let $H'$ be generated
by the unused positions. Lemma~\ref{chart-macro} gives
$|H'|\ge n^r/I$ as long as the unused length exceeds
$(r-1+1/2)n$. Every generator of $K$ is a difference of sums in
$H'$, so $K\le H'$. The preceding factorization yields
$|H'/K|\le n^k$, hence $|K|\ge n^{r-k}/I$.

At the start of each test take $H$ to be the subgroup $K$ obtained
before any endpoints have been used in that test. All exchange
differences have $F$-value zero, so $H\le\ker F$. The factorization
just proved gives $P:\mathbb Z^k\to G/H$ prescribing every class
image as $Pa_j$. Put $\Lambda=\ker P$. Since $nG=0$ and
$F(Pa)=Ba\bmod n$ in the induced map on $G/H$, we have
$n\mathbb Z^k\subseteq\Lambda\subseteq nB^{-1}\mathbb Z^k$.
For every $u\in(n^{-1}\Lambda)^*$ the assignment
$Pa\mapsto u\cdot a\bmod n$ is well-defined. Extend it from the
image of $P$ to $G/H$ by the character extension above, then pull
it back along $G\to G/H$.
During a test the differences still available generate a subgroup
$K\le H$, since every exchange on unused positions was available at
the start of the test.

{Let $N_j$ be the current number of available terms in class $j$.}
Test whether there is a real vector $x$ with
$\epsilon N_j\le x_j\le(1-\epsilon)N_j$ and $Ax\in\Lambda$.
The image of this box under $A/n$ is exactly the homothetic copy
of $Z$ of ratio $1-2\epsilon$ about its center. Failure therefore
gives the required lattice $\Gamma=n^{-1}\Lambda$.
At $k=0$ the test succeeds. In a successful test write $\lambda=Ax$.
Each coordinate has the tested slack $\epsilon\tau n/2$.
Suppose disjoint exchanges use at most $\epsilon\tau n/4$ endpoints
in total. Apply \eqref{chart-eq:count-rounding} with $b=U\lambda$
and $D=1$ to obtain $u= x+O(1)$ with $Au=\lambda$. Choose the
first arm of every exchange. In each class the number already used
and the number forbidden in the other arms are each at most
$\epsilon\tau n/4$. Since $u_j$ and $N_j-u_j$ are both at least
$\epsilon\tau n/2-O(1)$, we can fill the selected arms to exactly
$u_j$ positions using unused terms. At least one filler remains
in the selection under every later switch. The initial group sum
belongs to $H$, since its image in $G/H$ is $P\lambda=0$.
If the optional sums of exchange differences covered $H$, switches
could cancel that initial sum and yield a nonempty zero. Thus no
pool, or union of disjoint pools within the same endpoint budget,
can have all of $H$ as its optional sums.
\smallskip\noindent\emph{Growth of the exchange sums.}
For $k=r$ we have $|H|\le C$. Start with $Y=\{0\}$ and consume an
available exchange if $Y\cup(Y+w)$ strictly enlarges $Y$, replacing
$Y$ by this union. Each step adds an element, so at most $|H|-1$
bounded exchanges use $O(1)\le\eta n$ endpoints. The count argument
forbids $Y=H$. At termination every unused difference stabilizes $Y$.
The stabilizer is proper, since $h+Y=Y$ and $0\in Y$ imply $h\in Y$, whereas $Y\ne H$.
Discard the used endpoints and retest: the new working subgroup is
generated by those unused differences and is therefore proper in $H$.

Now suppose $k<r$ and put $d=r-k$. To build a pool of optional exchange
differences, at each step use the current set of all available
differences, which generates the subgroup $K$. Take generators successively from the available differences, each
outside the subgroup generated so far, until that subgroup is $K$.
Their relative orders $t_i$ divide $n$ and have product $|K|$.
Sums with digits $0,\ldots,t_i-1$ are distinct: reducing an equality
modulo the subgroup before the last generator forces equality of its
digits, and descending removes every generator in turn.
For $2\le L'\le n$, $\min(t_i,L')\ge t_i^{\log L'/\log n}$, by
considering $t_i\le L'$ and $t_i>L'$ separately. The product of
these upper digit lengths is at least
$(n^d/I)^{\log L'/\log n}\ge (L')^d/I$.
If it exceeds $(L')^d$, stop at the first factor which would exceed
that threshold, replace it by the integer part of the threshold
divided by the preceding product, and set later lengths to one.
Since $\lfloor a\rfloor\ge a/2$ for $a\ge1$, the resulting
product is at least $(L')^d/(2I)$ and at most $(L')^d$.
For each retained digit length $s>1$ use
$(s-1)/\log s\le (L'-1)/\log L'$, which follows by
differentiating $(s-1)/\log s$. Summing and taking the logarithm
of the product bound gives total digit cost at most $d(L'-1)$.

If the current optional sumset $Y$ has size $z$, take
$L'=\lceil(4Iz)^{1/d}\rceil$. There are at least $2z$ distinct
digit translations, and $\sum_{h\in H}|Y\cap(Y+h)|=z^2$
by counting ordered pairs from $Y$. One translation thus loses
at least $z/2$ elements of overlap. For sums of translations
the loss is subadditive, since twice the loss equals the
symmetric-difference distance. The digit cost bound gives one
available difference with gain at least $z/(2d(L'-1))\ge c_I
z^{1-1/d}$. Consume only that one exchange. Writing $z'$ for
the new size, $z\le z'\le2z$, and the mean value theorem gives
$(z')^{1/d}-z^{1/d}\ge(z'-z)/(d(2z)^{1-1/d})\ge c'_I$.
Hence $O_I(\rho^{1/d}n+1)$ exchanges suffice to reach
$z\ge\rho n^d$. Choose $\rho$ small enough that $L'\le n$
throughout this stage and the endpoints used are at most
$\eta n/2$ in each pool.
To grow this pool further, choose an available difference giving
relative gain at least $M/n$ at each step. Starting from density
at least $\rho n^d/|H|\ge\rho/C$ for a fixed $C$, after $b$ such
steps its density is at least $(\rho/C)(1+M/n)^b$. For a prescribed
small positive fraction $b/n$, take $M$ large enough that this
exceeds $1/2$ for large $n$. All steps fit the chosen endpoint
budget after allowing for the fixed maximum size of an exchange.
Repeat with fresh positions for a second pool. Two subsets of $H$
of size greater than $|H|/2$ have sum $H$, since $Y_1$ meets
$h-Y_2$ for every $h\in H$. {By the count test above,} one of
the two constructions must stop at a nonempty set $Y$ of density
$\theta\le1/2$, with $|(Y+w)\setminus Y|<M|Y|/n$ for every
remaining exchange difference $w$.

\smallskip\noindent\emph{From slow growth to a character.}
Use Haar probability on $H$ and put $f=1_Y$. {For a character $\chi$ written in circle notation, let $\widehat f(\chi)=|H|^{-1}\sum_{x\in H}f(x)\overline{\chi(x)}$.} The nontrivial Fourier
coefficients have squared moduli summing to
$\|f\|_2^2-|\widehat f(1)|^2=\theta-\theta^2$. Normalize these
squared moduli to a probability on the nontrivial characters.
For an integer $j\ge1$, the triangle inequality for symmetric
differences and the bound at $w$ give
$$
  |Y\mathbin\triangle(Y+jw)|\le j|Y\mathbin\triangle(Y+w)|
                      \le 2jM|Y|/n .
$$
Parseval applied to the difference of the two indicators expresses
the left side divided by $|H|$ as
$2\sum_\chi |\widehat f(\chi)|^2(1-\Re\chi(jw))$.
Dividing by $2\theta(1-\theta)$ now gives
$\mathbb E_\chi(1-\Re\chi(jw))\le 2jM/n$.

Write $\chi(w)=\exp(2\pi i t/n)$ with $|t|\le n/2$, and take
$J=\lfloor n/T\rfloor$ for a fixed $T$. When $|t|>4T$,
the geometric-sum formula and
$|\sin(\pi t/n)|\ge2|t|/n$ give
$$
 \left|J^{-1}\sum_{j=1}^{J}\chi(jw)\right|
       \le\frac{1}{J|\sin(\pi t/n)|}
       \le\frac{n}{2J|t|}<\frac14
$$
for large $n$. Averaging the preceding Parseval bound over
$1\le j\le J$ gives at most $M(J+1)/n$. Thus the probability of
$|t|>4T$ is at most $4M(J+1)/(3n)$. Average once more over any
chosen distribution of available exchanges. For $T$ sufficiently
large, some nontrivial character has $|t|\le4T$ on all but an
arbitrarily small fraction of that distribution.
\smallskip\noindent\emph{Augmentation or descent.}
Give the finitely many exchange types equal probability. Within each
type choose distinct positions uniformly in their designated classes.
In sampling its bounded number of slots independently the probability
of a repeated position is $O(1/n)$, since each class has at least
$\tau n/2$ terms. Thus making the overall exceptional probability
arbitrarily small also makes it arbitrarily small in each type,
with or without replacement.

Extend the character to $\ell$ on $G$. Call an exchange good when
its $\ell$-residue has circular distance at most $4T$ from zero.
Write $E=4T$ if any good exchange has nonzero residue, and $E=0$
otherwise. In either case all good residues have distance at most $E$.
In each class, averaging over the first term of a two-term exchange
gives a value $b_j$ within distance $E$ of all but an arbitrarily
small fraction of the terms. For each relation $w_j$, sample the
positions of its two arms. The sum of the failure probabilities
for its finitely many individual bounds and for goodness of the
exchange is less than one for sufficiently large $T$. Thus one
sample satisfies them all, and the triangle inequality gives
${\lVert b w_j\rVert_n}\le E(1+\|w_j\|_1)$, where $b=(b_j)_j$ is a row vector.
Discard the terms outside the individual bounds, within the single
centering budget of $\eta n$.
If $E=4T$, choose a good nonzero residue $t$. Then $0<{\lVert t\rVert_n}\le4T$,
so the order of $\ell|_H$ is at least
$n/\gcd(n,t)\ge n/(4T)$. Put $D=\det B$ and
$f=D\ell-(bU)\operatorname{adj}(B)F$, with $f=\ell$ at $k=0$.
For class $j$, $F=Ba_j$ and
$(bU)\operatorname{adj}(B)Ba_j=D\,bUa_j$. The center of $f$ on
this class is consequently $D(b_j-bUa_j)=D\,b w_j$, which has
bounded residue by the preceding bound. Each kept value of $f$
differs from that center by at most $4|D|T$. On $H$ one has
$f=D\ell$, so $|f(H)|\ge n/(4|D|T)$, taking $D=1$ at $k=0$.
The image of $(F,f)$ projects onto $\operatorname{im}F$ with a
kernel containing $\{0\}\times f(H)$. Its index in $C_n^{k+1}$
is therefore bounded. Refine the retained classes by their finitely
many $f$ values, and start the next chart.

If $E=0$, all retained values equal their class centers modulo $n$
and $b w_j=0\bmod n$ for every relation.
Every surviving exchange therefore has zero $\ell$-value: for
a same-class pair its two values agree, and for the arms of $w_j$
their difference is $b w_j=0$. Retesting therefore gives a new working
subgroup contained in $H\cap\ker\ell<H$, since $\ell$ is nontrivial on $H$.

Write $H_t$ for the working subgroup recomputed after $t$ proper
descents. For $t\le L$ the remaining classes and length
satisfy the invariants above. Their available differences still
generate $K\le H_t$ with $|K|\ge n^{r-k}/I$. Therefore
$$
       n^{r-k}/I\le |H_t|\le 2^{-t}|\ker F|\le 2^{-t} C n^{r-k}.
$$
This gives $t\le\log_2(CI)<L$. Thus at most $L$ tests are needed
at a given chart, and the budgets allocated above suffice without
assuming the eventual outcome of any test.
In particular, at $k=r$ a successful test can only cause a
descent, so a count test must fail within this budget. At any
dimension the resulting lattice $n^{-1}\Lambda$ lies between
$\mathbb Z^k$ and $B^{-1}\mathbb Z^k$. The latter quotient by
$\mathbb Z^k$ is finite, giving only finitely many choices.
At each of at most $r$ augmentations the bounded new coordinate
values allow only finitely many new column lists. Taking the
largest constants over the finitely many possibilities makes
all stage bounds uniform at fixed $\epsilon$. The last test
fails in a positive dimension because the zero-dimensional
test always succeeds.
\end{proof}

We now let $\epsilon\to0$. The chart constants may grow, so we
normalize $Z$ before choosing the lattice directions to retain.
\begin{lemma}\label{chart-f:compact}
Along a subsequence of any family of zero-sumfree sequences {$(S_n)$ {in $G=C_n^r\oplus E_n$} with
$|S_n|/n\to r$}, there are $1\le q\le r$ and homomorphisms
$F_n:G\to C_n^q$ of uniformly bounded image index. On all but
$o(n)$ positions their outputs have integer representatives in $\mathbb Z^q$
whose normalized counting measures converge pointwise to weights
$\alpha_v\ge0$. They have uniformly bounded first moment, and
$$
 \sum_v\alpha_v=r,\qquad \sum_v\alpha_v\|v\|<\infty,\qquad
 \operatorname{span}\{v:\alpha_v>0\}=\mathbb R^q .
$$
The body $\sum_v[0,\alpha_v]v$ is hollow in $\mathbb Z^q$.
\end{lemma}
\begin{proof}
Let tolerances decrease to zero, at each fixed tolerance taking an
exponent sufficiently large for the estimates below. At tolerance
$\delta_j$ let $c_j>0$ be the uniform lower bound for $N_i/n$ from
Lemma~\ref{chart-f:chart}. {For each of the finitely many column lists and lattices at this tolerance, choose, for every nonspanning sublist, a nonzero covector in the original dual lattice $\Gamma^*$ that vanishes on the sublist. Let $B_j$ bound the absolute evaluations of these covectors on their full lists.} Choose $n_j$ also so that $n_j c_j\ge j$ {and $n_j>I(1/2)B_j$}. Pass to a
subsequence of fixed terminal dimension $k$. Use as coordinate basis
$k$ vectors of maximal absolute determinant in the centered $Z$.
In this basis the centered body contains the unit crosspolytope by
convexity and symmetry, and is in the unit cube by Cramer's rule.
Transform the original body and lattice by this same linear map.
As $0\in Z$, the bodies stay bounded. Moreover the sum of their
coordinate widths is $\sum_i(N_i/n)\|a_i\|_1\le2k$.
Consequently all but at most $n/4$ retained positions have norm
at most $8k$. If their columns failed to span, choose a normal in
$\Gamma^*$ killing those columns but not all columns. This is possible
because the columns are lattice vectors and the full list spans.
Some bounded nonzero evaluation $m$ of this normal would give an
actual character of order at least $n/|m|$, whereas the number of
annihilated terms exceeds $(r-1+1/2)n$. Here the integer normal
and its evaluations can be bounded at each fixed tolerance by choosing
one normal for each of the finitely many original column subsets. {With this choice, $n_j/|m|\ge n_j/B_j>I(1/2)$, contrary to the character bound.}
We thus have $k$ independent lattice vectors of bounded length.
Write $Z_j,\Gamma_j,\delta_j$ for the normalized body, lattice and
tolerance at step $j$. In each $\Gamma_j$,
choose a shortest nonzero vector and project along it. The vector is
primitive, since division by any common factor would shorten it. {Identify each real quotient with the orthogonal complement of this vector and use the Euclidean distance there.}
It extends to a lattice basis, whose remaining projections generate
the quotient lattice freely. Repeat in the quotient, lifting each
chosen vector back and subtracting earlier basis vectors in reverse
order to put earlier Gram--Schmidt coefficients in $[-1/2,1/2]$.
This gives a basis $b_{j,1},\ldots,b_{j,k}$ with Gram--Schmidt
lengths $d_1,\ldots,d_k$. Each $d_i$ is uniformly bounded above:
one of the $k$ bounded independent vectors projects nontrivially
past the first $i-1$ vectors. At that step the projection of
$b_{j,i+1}$ has squared length between $d_i^2$, by shortestness,
and $d_i^2/4+d_{i+1}^2$, by reduction. Hence
$d_{i+1}\ge\sqrt3\,d_i/2$. Pass to limits of these lengths.
Those tending to zero form an initial segment $1,\ldots,k-q$.
Size reduction makes the initial basis vectors tend to zero as
well. Rounding the coefficients in this basis covers their real
span to within half the sum of their lengths, hence to within $o(1)$.
If $q=0$, rounding the center would put a lattice point inside the
fixed-radius ball contained in the contracted body. Thus $q\ge1$.
Project orthogonally off that first span and express the projection
in the basis of the last $q$ projected lattice vectors. This linear
map $T_j$ sends $\Gamma_j$ onto $\mathbb Z^q$ and has bounded
operator norm. On the orthogonal complement its inverse also has
bounded norm. This follows from the two-sided Gram--Schmidt
bounds and size reduction. The width estimate now controls every column individually:
$(N_i/n_j)\|a_i\|_1\le2k$ gives $\|a_i\|_1\le2k/c_j$ in the
normalized coordinates. Hence for each fixed integral row $\ell$,
$$
 \max_i\frac{|\ell T_j a_i|}{n_j}
 \le\frac{2k\|\ell\|\,\|T_j\|}{n_j c_j}
 \longrightarrow0 .
$$

Take a Hausdorff limit $K$ of $T_jZ_j$. It has nonempty interior
since the transformed bodies contained fixed balls. Its interior
misses $\mathbb Z^q$. Indeed if an integer $u$ lay in the interior,
it would belong, for large $j$, to the image of a fixed smaller
concentric homothety of $Z_j$, by support-function convergence.
Choose a real lift there and any lattice lift of $u$. Their
difference is in the short span, where the short lattice has
covering radius $o(1)$. Thus a lattice point in the same fiber
is within $o(1)$ of the real lift. The fixed homothety leaves
a positive Euclidean margin inside the forbidden contracted body,
a contradiction.

{Pull the rows of $T_j$ back through the normalizing linear map from the original lattice to $\Gamma_j$. Since each row is integral on $\Gamma_j$, the resulting covectors belong to the original dual lattice.} By Lemma~\ref{chart-f:chart}, {they are realized on the retained classes by actual characters.} {Let $F_{n_j}:G\to C_{n_j}^q$ have these characters as its coordinates.} Take the retained integer columns with
normalization $h_j=n_j$. Their normalized count tends to $r$.
The width bound $2k$ and boundedness of $T_j$ give the required
first-moment bound in Lemma~\ref{chart-l:column-limits}. Write $\alpha_v$
for its limiting weights. The approximating subsums lie in $T_jZ_j$
after division by $n_j$. Hausdorff convergence places their limits
in $K$, proving hollowness of the closed sum of the limiting segments.
{Call a column $v$ occupied if $\alpha_v>0$.} A proper span of the occupied columns would admit a fixed nonzero
integer normal $\ell$. By tightness, ${\ell F_{n_j}}$ would kill
$rn_j-o(n_j)$ terms. Yet the projected columns at each $j$ span
$\mathbb R^q$, so some evaluation $m$ of $\ell$ is nonzero.
By the width estimate, $0<|m|/n_j\to0$, even if the
chosen column varies with $j$. Its residue has order
$n_j/\gcd(n_j,m)\ge n_j/|m|\to\infty$. But the character bound
preceding Lemma~\ref{chart-f:chart} gives $|({\ell F_{n_j}})(G)|\le I(1/2)$
from the number of annihilated terms. This is a contradiction.
Hence the limiting columns span.
In particular, choose $q$ fixed independent
occupied columns with square column matrix $D_0$. Each occurs
for all large $j$. Their integral span has index $|\det D_0|$
in $\mathbb Z^q$ and its reduction modulo $n_j$ is in the image
of the actual coordinates. The image therefore has index at most
$|\det D_0|$, uniformly in $j$.
\end{proof}

{\vspace{0pt}
\begin{proof}[Proof of Theorem~\ref{chart-main}]
Return to the notation $G_n=C_n^d\oplus E_n$ of the theorem. By Lemma~\ref{chart-f:compact}, we may pass to a subsequence to obtain some $k\ge1$
and homomorphisms $F_n:G_n\to C_n^k$ with bounded image index.
On all but $o(n)$ positions the outputs have integer representatives.
Their normalized counting measures are tight and converge pointwise
to a measure $\mu$ of mass $d$ on $\Z^k$ with finite first moment.
The occupied columns, those with $\mu(v)>0$, span $\R^k$ and the compact convex set
$$
 Z=\sum_{v\in\Z^k}[0,\mu(v)]v
$$
has interior disjoint from $\Z^k$. The sum denotes all sums of
independently chosen coefficients, which converge absolutely by the
first-moment bound.

For $0\le s\le d$ and any homomorphism $L:G_n\to C_n^s$ of bounded image index, restriction
to $C_n^d$ increases that index by at most $|E_n|$. Smith normal
form of an integral lift of this restriction gives $d-s$ free
$C_n$ coordinates in the kernel. The other kernel factors have
orders $\gcd(n,a_i)$, where the $a_i$ are the diagonal entries in Smith normal form. Their product equals the restriction's image
index. The $d-s$ free coordinates constitute a direct summand of
$G_n$ contained in $\ker L$. Its complement in $\ker L$ has
bounded order, by $|\ker L|=|G_n|/|\operatorname{im} L|$.
This also covers $s=0$ and $s=d$.

For a subspace $W$, put $\mu(W)=\sum_{v\in W\cap\Z^k}\mu(v)$. If $W$ is rational of dimension $u$, project along it
using the surjective lattice map with kernel $W\cap\Z^k$.
This primitive quotient is surjective modulo every $n$ (extend a
basis of the saturated kernel to a lattice basis).
Its composition with $F_n$ thus has bounded image index.
All the positions represented by a vector in $W$ lie in its kernel.
Testing first finitely many occupied columns, applying Lemma~\ref{chart-macro}
to that kernel and then increasing the finite set gives
\begin{equation}\label{chart-cap}
                     \mu(W)\le d-k+u.
\end{equation}
The same holds for real $W$: its contained lattice columns span a
rational subspace of dimension at most $u$.

\medskip\noindent\textit{A terminal quotient and its complementary subspaces.}\quad
Among the rational quotient projections of the hollow body, including
the identity, keep a hollow one of least positive dimension. Write $k$ for this dimension, identify its quotient lattice with $\Z^k$, and push $\mu$ forward. The new
characters still have bounded image index. Tightness and finite first
moment persist, the columns span, and \eqref{chart-cap} holds by the same
kernel argument. Every proper {rational} quotient of this new space has a
lattice point in the relative interior of the projected body; this
includes the zero-dimensional quotient.

We first prove $\mu(0)=d-k$. Say a rational subspace $W$ is tight
if equality holds in \eqref{chart-cap}. The full space is tight. Choose
a tight $W$ of smallest dimension $u$. If $u>0$, every proper
subspace $V$ of $W$ has $\mu(V)<d-k+\dim V$: for real $V$ test the
rational span of its contained columns. Upper semicontinuity on each
Grassmannian gives uniform slack by compactness. In $W$, the mass
outside each proper $V$ therefore exceeds $u-\dim V$ by a uniform
positive amount. Retain finitely many columns in $W$ with total
omitted mass smaller than this slack. These columns span $W$ and
their positive open segment sum covers $W/(W\cap\Z^k)$ by the
finite covering lemma (Lemma~\ref{cov-cover}).
Its hypotheses are exactly the just-verified outside-mass inequalities
on deficient sublists, with the columns written in any basis of
$W\cap\Z^k$. Zero columns do not contribute to outside mass.

Let $C'$ be the closed segment sum of the columns outside $W$,
with real span $U$. Projection $\pi$ along $W$ maps $C'$ onto
the whole projected body. By terminality pick an interior projected
lattice point and lift it to $x\in\operatorname{relint}_U C'$.
For this lift, use the segment from the relative interior center of
$C'$ towards a lift of a projected point slightly beyond the
desired point; positive weight on the center gives relative interior.
If the desired point is the projected center use the center directly.
Let $z\in\Z^k$ be any lattice lift of the same point. Since
$z-x\in W$, open coverage gives an interior point of the segment
sum in $W$ equal to $z-x$ modulo $W\cap\Z^k$.
Adding $x$ gives a lattice point in the interior of $Z$: the two
relative neighborhoods have sum a neighborhood in $U+W=\R^k$.
This contradiction forces $u=0$, proving the mass at zero.

Discard the zero column in the measure. The remaining mass is $k$
and every subspace has mass at most its dimension. Choose a
nonzero rational tight subspace $W$ minimal by inclusion, now
using this remaining mass. Its columns span $W$, have mass $\dim W$,
and satisfy the strict subspace inequalities there. Define $C',U,x,z$
as above. For every sufficiently small $b\in U\cap W$ the point
$x+b$ is also in $\operatorname{relint}_U C'$. The target
$z-x-b$ modulo $W\cap\Z^k$ must be missing from the open segment
sum in $W$, or the preceding addition would violate hollowness.
Lemma~\ref{chart-isolated} gives finitely many missing targets.
A positive-dimensional neighborhood in $U\cap W$ takes infinitely
many values modulo a discrete lattice. Hence $U\cap W=0$.
Thus $U$ is a rational complement, its remaining mass equals
$k-\dim W=\dim U$, and the family in $W$ has the asserted
nearest-point description.

If $U$ is nonzero, choose a minimal nonzero tight subspace inside
$U$. It is also minimal among tight subspaces contained in it in
the full space. Apply the same argument with its full complement
of occupied columns, including the already extracted families.
Its intersection with the span of all those other columns is zero.
Iterating decreases the remaining dimension and proves the direct
decomposition and finiteness. Finite limiting support of full mass
implies that all but $o(n)$ actual terms have the listed outputs.
Therefore precisely $(d-k)n+o(n)$ terms are in the actual kernel.
For $k=1$ the hollow interval contains $0$ and has nonzero column
mass one. Its length is $\sum|\!v\!|\mu(v)\ge1$ and at most one by
hollowness. Every nonzero integer of positive limiting mass is thus $\pm1$.
Both signs cannot have positive mass since $0$ would be interior.
The value of positive mass occurs for all sufficiently large $n$, so $F_n$
is then surjective.
\end{proof}

The kernel subsequence in part 1 of Theorem~\ref{chart-main} again satisfies the hypotheses if
$d-k>0$. Iteration terminates in at most $d$ applications. This
controls the projections of the successive groups of positions. The
components of an upper group in earlier kernels can still vary.
\par}

\section{Prime-power cancellation and ordinary flag stability}\label{sec:stability}

We shall prove in Theorem~\ref{stability} that a zero-sumfree family of length
$rn+o(n)$ admits, on a subsequence, an ordinary flag with $n+o(n)$
positions in each row. We first establish exact cancellation when the
total deficit is smaller than a prime-power divisor, then lift a strict
selection across a kernel and lift the resulting flag to the original
modulus.

\subsection{Exact cancellation below a prime-power deficit}\label{sec:small}

Write $C_n=\Z/n\Z$. A sequence is a list of indexed positions. For a basis $e_1,\ldots,e_r$ of $C_n^r$, its $j$th ordinary flag row is
$$
 e_j+\langle e_1,\ldots,e_{j-1}\rangle .
$$
We show directly that a bounded deficit in the stocks of such rows cannot sustain a critical zero-sumfree sequence. The lower coordinates of each term are arbitrary and may vary from position to position.

\begin{theorem}\label{small-main}
Let $n\ge2$, $r\ge1$, and $Q=p^a\mid n$ with $p$ prime and $a\ge1$. Let $S$ have length $T_r(n)=r(n-1)+1$. Suppose there are $n-1-d_j$ indexed terms chosen in its $j$th ordinary flag row for $j=1,\ldots,r$, where $0\le d_j\le n-1$. Put $d=\sum_{j=1}^r d_j$. If $d<Q$, then $S$ contains a nonempty zero-sum subsequence.
\end{theorem}
There are exactly $d+1$ positions outside the chosen rows; there are no restrictions on those positions. Excess terms in a row may be included among them. The proof pads the rows to size $n-1$ with artificial terms and gives weight zero to every correction using any of them. A combinatorial averaging functional controls the degree of this weight even when an ordinary uniform average would have a denominator divisible by $p$.

\medskip\noindent\textit{Averaging on a slice of the Boolean cube.}\quad
Write $[L]=\{1,\ldots,L\}$, with $[0]=\varnothing$, and let $\Omega(L,k)$ be the set of its $k$-element subsets. For $A\subseteq[L]$ write $x_A$ for the function $x_A(U)=1$ if $A\subseteq U$ and $0$ otherwise. These are the squarefree monomials on the indicator vectors of the subsets. Call $A=\{a_1<\cdots<a_j\}$ admissible if $a_i\ge2i$ for all $i$; the empty set is admissible. Every subset of an admissible set is admissible.

\begin{lemma}\label{small-basis}
The functions $x_A$ for admissible $A$ with $|A|\le\min(k,L-k)$ form an integral basis of the integer-valued functions on $\Omega(L,k)$.
\end{lemma}
\begin{proof}
Proceed by induction on $L$. When $k=0$ the single constant function is a basis, including $L=0$. At a given $L$ first take $1\le k\le L/2$. On the sets not containing $L$, the functions indexed by $A$ not containing $L$ restrict exactly to the inductive basis on $\Omega(L-1,k)$. Indeed, for $L>2k$ its size cutoff is $k$, while for $L=2k$ it is $k-1$, and an admissible $k$-set cannot fit in $[2k-1]$. Each remaining admissible index has the form $A=B\cup\{L\}$ where $B$ is admissible with $|B|\le k-1$; the last condition for admissibility holds because $L\ge2k$. Its function vanishes on every $U$ not containing $L$. On the sets $U$ containing $L$ it is exactly $x_B(U\setminus\{L\})$, giving the inductive basis on $\Omega(L-1,k-1)$. One can thus represent any integer function uniquely, first matching its values on sets without $L$ and then its remaining values on sets containing $L$, with integral coefficients. For $k>L/2$ pass to complements and use the already proved case $L-k$. Under complementation $x_A$ becomes $\prod_{i\in A}(1-x_{\{i\}})=\sum_{B\subseteq A}(-1)^{|B|}x_B$. The index sets are closed under subsets and the transformation is triangular by size with diagonal entries $\pm1$. This proves the assertion.
\end{proof}

\begin{lemma}\label{small-moment}
Suppose $L=n-1$, $Q=p^a\mid n$, $0\le D<Q$, and $0\le k\le L$. On the $\F_p$-valued functions on $\Omega(L,k)$ that are restrictions of multilinear polynomials of degree at most $D$, there is a linear functional $M_k$ such that for every $I\subseteq[L]$ with $t=|I|\le D$,
$$
 M_k(x_I)=\bin{k}{t}\big/\bin{L}{t}\quad\text{in }\F_p.
$$
In particular this sends the constant one to one and any function vanishing identically on the slice to zero.
\end{lemma}
\begin{proof}
For $t<Q$, the coefficient of $z^t$ in $(1+z)^L$ equals $(-1)^t$ in $\F_p$. Indeed $(1+z)^n=(1+z^Q)^{n/Q}$ is one to degrees below $Q$, and $(1+z)^L=(1+z)^n/(1+z)$ in formal power series. Thus every denominator in the formula is a unit. If $\min(k,L-k)<D$, then $\bin{L}{k}=\bin{L}{\min(k,L-k)}$ is itself a unit, so ordinary averaging over the slice works. The desired ratio follows by counting pairs $(I,U)$ with $I\subseteq U$, or as zero if $t>k$.

Assume $D\le\min(k,L-k)$. On the basis from Lemma~\ref{small-basis}, prescribe the stated moments for all basis indices of size at most $D$ and any values on the larger ones. Reduction of the integral basis modulo $p$ is a basis, so this defines a functional. We must verify its values on every monomial, including those not in the basis. Fix $0\le t\le D$. By the same lemma on the $t$-subset slice, the formal coefficient vectors of the sums
$$
 P_{A,t}=\sum_{\substack{I\supseteq A\\ |I|=t}}x_I,
 \qquad A\ \hbox{admissible},\quad |A|=j\le t ,
$$
are an integral basis for all coefficient vectors on $t$-subsets: here $t\le L/2$. On $\Omega(L,k)$ the displayed polynomial equals $\bin{k-j}{t-j}x_A$ since, once $A\subseteq U$, it counts the extensions in $U$. Consequently the chosen functional sends $P_{A,t}$ to $\bin{k-j}{t-j}\bin{k}{j}/\bin{L}{j}$. The proposed moments term by term give $\bin{L-j}{t-j}\bin{k}{t}/\bin{L}{t}$, the same rational number by the factorial formulas. Both rational expressions have denominators prime to $p$ and hence have equal reductions. Equality on the coefficient basis proves the desired formula on every individual $x_I$.
\end{proof}

\medskip\noindent\textit{Correcting the heights without using an artificial term.}\quad
In variables $z=(z_1,\ldots,z_s)\in\Z^s$, let $\mathcal P_b$ be the $\F_p$-span of the functions $\prod_i\bin{z_i}{u_i}\bmod p$ with nonnegative $u_i$ and $\sum_i u_i\le b$. Set $\mathcal P_b=0$ for $b<0$. Binomial coefficients with negative upper arguments are taken in their usual integer-valued polynomial sense. If $b<Q$, these functions have period $Q$ in each variable: in the generating function $(1+v)^{z_i}$ a shift by $Q$ multiplies by $1+v^Q$ modulo $p$, unchanged through degree $Q-1$. Thus they also define functions on $C_n^s$. For any integral shift $w$, $\Delta_w f(z)=f(z+w)-f(z)$ takes $\mathcal P_b$ into $\mathcal P_{b-1}$. Indeed the expansion of each translated binomial by Vandermonde's identity has just one term of its original multi-index, and it cancels. In one variable $\bin{-z}{t}=(-1)^t\bin{z+t-1}{t}$ for $t\ge1$, so $\bin{-z}{t}\in\mathcal P_t$ by Vandermonde. For $t=0$ the function is one.

\begin{proof}[Proof of Theorem~\ref{small-main}]
Set $L=n-1$. Add $d_j$ distinguishable artificial copies of $e_j$ to the chosen terms in row $j$. This yields $L$ labeled slots per row, all with $j$th coordinate one and all higher coordinates zero. Let $D_j=\sum_{i\le j}d_i$. From a state $z\in C_n^j$, full correction beginning at row $j$ chooses exactly $k=[-z_j]_n$ slots in that row, where brackets denote the representative in $[0,n-1]$. This cancels the $j$th coordinate. {Write $z_{<j}=(z_1,\ldots,z_{j-1})$.} With $b_i$ the indicator of the chosen slots and $w_i$ their vectors of lower coordinates, the state for the next row is $z_{<j}+\sum_i b_i w_i$. A correction is legal when none of its selected slots is artificial.

{Let $J_j$ index the $d_j$ artificial slots in row $j$.}
Starting with $P_0=1$ on the trivial state, define $P_j(z)$ by applying $M_k$ of Lemma~\ref{small-moment} at row $j$ to
$$
 \prod_{i\in J_j}(1-b_i)\ P_{j-1}\left(z_{<j}+\sum_i b_iw_i\right),
$$
We verify inductively that this definition is meaningful, that $P_j\in\mathcal P_{D_j}$, that $P_j(0)=1$, and that $P_j(z)=0$ whenever $z$ admits no legal correction down through row one.

{Put $b_I=\prod_{i\in I}b_i$ for $I\subseteq[L]$.} By translation and expansion on Boolean variables,
$$
 P_{j-1}\left(z_{<j}+\sum_i b_iw_i\right)=
 \sum_{I\subseteq[L]} b_I\left(\prod_{i\in I}\Delta_{w_i}\right)P_{j-1}(z_{<j}).
$$
This follows by multiplying the commuting translation operators $1+b_i\Delta_{w_i}$. The coefficient at $b_I$ is in $\mathcal P_{D_{j-1}-|I|}$. Multiplication by the artificial-slot product and reduction by $b_i^2=b_i$ gives at $b_T$ a coefficient in $\mathcal P_{D_j-|T|}$, since $|T|\le |I|+d_j$ for every contributing $I$. Thus only $|T|\le D_j<Q$ can occur and Lemma~\ref{small-moment} applies. It replaces $b_T$ by $\bin{k}{|T|}/\bin{L}{|T|}$. Periodicity replaces $k$ by $-z_j$ in the numerator modulo $p$. The resulting factor in the new variable has degree at most $|T|$ in the binomial basis. Multiplying by the coefficient in the other variables gives total degree at most $D_j$, as required. At $z=0$, $k=0$ and the unique selection is empty, giving $P_j(0)=P_{j-1}(0)=1$. If $z$ admits no legal correction, each $k$-selection either uses an artificial slot or reaches a state with no legal correction. The function being averaged is then zero at every point of the slice, so $P_j(z)=0$ by the lemma.

Let $u_1,\ldots,u_{d+1}$ be the actual terms outside the chosen rows. Under zero-sumfreeness, a nonempty selection of them cannot be legally corrected. Hence $P_r(\sum_{i\in I}u_i)=0$ for every nonempty $I\subseteq[d+1]$, while $P_r(0)=1$. But
$$
 0=(\Delta_{u_1}\cdots\Delta_{u_{d+1}}P_r)(0)
 =\sum_{I\subseteq[d+1]}(-1)^{d+1-|I|}
 P_r\left(\sum_{i\in I}u_i\right)=(-1)^{d+1}\quad\text{in }\F_p.
$$
The initial zero follows from $P_r\in\mathcal P_d$. This contradiction proves the theorem, also when $d=0$.
\end{proof}

\subsection{Lifting a strict selection across the kernel}\label{sec:slice}
Let $q=p^a\to\infty$ be prime powers, allowing both $p$ and $a$ to
vary. Write $C_q=\Z/q\Z$. In $C_q^r$ let $U\simeq C_q^d$ be a
direct summand, with a chosen complementary quotient basis. Let $T$
be a sequence in $U$ with $|T|=d(q-1)-o(q)$. Outside this indexed
sequence there are disjoint classes of positions indexed by a fixed
finite set. Their projections to $C_q^r/U$ are fixed integral
columns $v_i\bmod q$ in the quotient basis, and class $i$ contains
$\lambda_i q+o(q)$ positions, $\lambda_i>0$. Within a class the
individual $U$-coordinates can vary arbitrarily. The bases, positions
and their $U$-coordinates may all depend on $q$. {For a sublist $J$, write $V_J$ for the matrix with columns $v_i$, $i\in J$.}

\begin{theorem}\label{slice-main}
Suppose on a nonempty sublist $J$ there are fixed real
$0<t_i<\lambda_i$ with
$$
       V_Jt=\sum_{i\in J}t_i v_i\in\sum_{i\in J}\Z v_i .
$$
Then, for every sufficiently large $q$, $T$ together with the
positions in the classes of $J$ admits a nonempty zero-sum
subsequence. If $T$ is zero-sumfree, one can use a positive
fraction of the positions in each class of $J$, leaving a positive
fraction unselected there.
\end{theorem}

{We prove Theorem~\ref{slice-main} after establishing the \hyperref[lem:prescribed-counts]{lemma on correction at prescribed counts}.}

\begin{lemma}\label{slice-residual}
Suppose $T$ is zero-sumfree and let $e=d(q-1)-|T|$. Then $e\ge0$.
For $e<q$ there is a function $F:C_q^d\to\F_p$ in the span of
$$
          \prod_{j=1}^d\binom{z_j+\alpha_j}{\alpha_j}\bmod p,
          \qquad \alpha_j\ge0,\quad\sum_j\alpha_j\le e ,
$$
with $F(0)=1$. Whenever $F(z)\ne0$, a subsequence of $T$ has sum
$z$, with the empty subsequence allowed.
\end{lemma}
\begin{proof}
In $\F_p[C_q^d]$ write $Y_j$ for the coordinate basis elements and
$X_j=1-Y_j$. Then $X_j^q=0$ and the $X_j$ generate this algebra:
$Y_j=1-X_j$ generate all group elements. The $q^d$ monomials of
individual exponents below $q$ form a basis by dimension. For an
element $g$ write $Y^g$ for the corresponding group element. Each
factor $1-Y^g$ is in the ideal generated by the $X_j$. {Writing $T=t_1\cdots t_s$ with $s=|T|$, the indexed product}
{$P=\prod_{i=1}^{s}(1-Y^{t_i})$} lies in its $|T|$th power. This product
has coefficient one at the identity since only the empty selection
has sum zero. Thus $P\ne0$, and the maximum possible monomial degree
$d(q-1)$ gives $e\ge0$. Each surviving term is a multiple of a
monomial $\prod_j X_j^{q-1-\alpha_j}$ with $\sum_j\alpha_j\le e$.
For the coefficient calculation first take $Y$ to be an
indeterminate in $\F_p[[Y]]$. The polynomial on the left has
degree below $q$, so its coefficients at degrees $0,\ldots,q-1$
are unchanged on then passing to the cyclic group algebra. We have
$$
        (1-Y)^{q-1-\alpha}=(1-Y^q)(1-Y)^{-1-\alpha}.
$$
Its coefficient at $Y^z$, for $0\le z<q$, is
$\binom{z+\alpha}{\alpha}$. For $\alpha<q$ this function of the
integer $z$ has period $q$ modulo $p$: multiplication of its
binomial generating function by $(1+x)^q=1+x^q$ preserves all
coefficients of degree below $q$. Defining $F(z)$ as the coefficient
of $Y^z$ in $P$ proves the claims.
\end{proof}

\begin{lemma}\label{slice-trim}
{Let $e=e(q)\ge0$ be integral with $e=o(q)$, and let the available stock be}
$N^*=\lambda q+o(q)$, $\lambda>0$. {For all sufficiently large $q$ one can retain
$N$ positions with $0\le N^*-N\le 2e$} such that all $\binom N t$ for
$0\le t\le e$ are units modulo $p$.
\end{lemma}
\begin{proof}
For $e=0$ use the entire stock. Otherwise let $L=p^j$ be the largest
power of $p$ at most $e$, and put $E=\lfloor e/L\rfloor\le p-1$.
Choose the largest $M\le N^*$ with $M\equiv-1\bmod L$.
If its digit $b=\lfloor M/L\rfloor\bmod p$ is at least $E$, use
$N=M$. Otherwise use $N=M-(b+1)L$. This changes digit $j$ to
$p-1$, preserves the lower digits $p-1$, and costs at most $EL\le e$
additional positions. The result is positive for large $q$ because
$N^*-N\le 2e$. Every digit of any $t\le e$ is dominated by the
corresponding digit of $N$, since $t<pL$. {Write $N=\sum_k N_k p^k$ and $t=\sum_k t_k p^k$ in base $p$.} Writing $(1+x)^N$ as
the product of $(1+x^{p^k})^{N_k}$ modulo $p$ gives
$\binom N t\equiv\prod_k\binom{N_k}{t_k}\ne0$.
\end{proof}

We use the following form of the Boolean slice functional from Lemmas~\ref{small-basis} and~\ref{small-moment}. For $0\le k\le N$, if
$\binom N t$ is a unit for every $0\le t\le e$, there is a linear
functional on the restrictions of squarefree polynomials of degree
at most $e$ to $k$-element subsets of $[N]$ which sends the monomial
on any $t$ distinct positions to
$\binom k t/\binom N t\bmod p$. The proof given above for $N=n-1$
uses the assumption on $N$ only to guarantee these units. Explicitly,
if $\min(k,N-k)<e$, then $\binom N k$ is also a unit and uniform
averaging gives the formula. Otherwise prescribe the moments on the
integral basis of the slice indexed by {admissible} subsets. To check every
degree-$t$ monomial, sum them over all $t$-sets extending one such
index $A$ of size $j\le t$. On the slice this sum equals
$\binom{k-j}{t-j}x_A$, with $x_A$ the monomial on $A$.
The prescribed and the desired moments on this sum agree because
$$
 \binom{k-j}{t-j}\frac{\binom k j}{\binom N j}
       =\binom{N-j}{t-j}\frac{\binom k t}{\binom N t}.
$$
The extension sums for {admissible} indices form the integral basis on
the $t$-subset slice, so they determine every degree-$t$ moment.

\begin{auxiliarylemma}[{Correction at prescribed counts}]
\leavevmode\phantomsection\label{lem:prescribed-counts}
{Let $T$ be zero-sumfree in $U=C_q^d$ with $0\le e=d(q-1)-|T|<q$.
In disjoint classes outside $T$, indexed by a finite list $J$, retain $N_i$ positions per class. Their
$U$-components may be arbitrary. Suppose $\binom{N_i}{t}$ is a unit modulo
$p$ for every $0\le t\le e$. There is an $\mathbb F_p$-linear
combination $A(k)$ of the functions $\prod_{i\in J}\binom{k_i}{s_i}$
with $\sum_i s_i\le e$ such that $A(0)=1$. If $0\le k_i\le N_i$ for all
$i$ and $A(k)\ne0$, one can choose exactly $k_i$ positions from each
class so that their total $U$-component is canceled by a subsequence
of $T$, possibly empty. In particular, for admissible counts this holds whenever each $k_i$
is divisible by a power $p^b>e$.}
\end{auxiliarylemma}
\begin{proof}
{Take the function $F$ supplied by Lemma~\ref{slice-residual}.}
Write $u_{ij}\in U$ for the $j$th retained lower component in class
$i$ of size $N_i$. For Boolean indicators $x_{ij}$ expand
$$
                      F\left(-\sum_{i,j}x_{ij}u_{ij}\right)
$$
in squarefree monomials over $\F_p$. Its total degree is at most
$e$. Indeed the coefficient at a set of $b$ distinct positions is
the mixed difference of $F$ in their $b$ negative directions at
zero. An integer translation lowers the total degree of an
integer-valued polynomial by at least one after subtracting the
original polynomial. Thus every mixed difference of order greater
than $e$ vanishes even before reduction modulo $p$, for each
binomial monomial in Lemma~\ref{slice-residual}.

Replace each Boolean monomial by the product over its classes
of $\binom{k_i}{s_i}/\binom{N_i}{s_i}$, where $s_i$ is the number
of positions of that class in the monomial. Let $A(k)$ denote the
resulting function of the counts. It is an $\F_p$-linear
combination of products of integer-valued binomial polynomials of
total degree at most $e$, and $A(0)=F(0)=1$. At any feasible
integer counts $0\le k_i\le N_i$, $A(k)$ is the value of the
tensor product of the slice functionals applied to the displayed
Boolean function. In particular $A(k)\ne0$ certifies an actual
choice with those counts for which $F$ is nonzero. This implication
uses only linearity on functions on the slices.
{The implication for $T$ follows from the defining property of $F$.
For the last assertion, each binomial of order at most $e$ has period
$p^b$ modulo $p$, since $(1+z)^{p^b}=1+z^{p^b}$ in characteristic $p$.
Thus $A(k)=A(0)$ for the indicated counts.}
\end{proof}

\begin{proof}[Proof of Theorem~\ref{slice-main}]
{If $T$ has a nonempty zero-sum subsequence the first conclusion
already holds. Assume henceforth that $T$ is zero-sumfree. By
Lemma~\ref{slice-residual} we have $0\le e=d(q-1)-|T|=o(q)$, so $e<q$
eventually. Trim each class of $J$ by Lemma~\ref{slice-trim} and let
$A$ be the function given by the {\hyperref[lem:prescribed-counts]{lemma on correction at prescribed counts}}.}
Write $V=V_J$ and fix $b\in\Z^J$ with $Vb=Vt$. Let
$u_1,\ldots,u_h$ be an integral basis of $\ker_{\Z}V$.
They span $\ker_{\R}V$, so $t-b=\sum_{\ell=1}^h\theta_\ell u_\ell$
with real $\theta_\ell$. Set
$$
        k=qb+\sum_{\ell=1}^h m_\ell u_\ell,\qquad
                  |m_\ell-q\theta_\ell|<\delta q .
$$
All these $k$ have quotient sum $0\bmod q$. Since $J$ and the
$u_\ell$ are fixed, $\delta>0$ can be chosen small enough that
every resulting $k_i$ and $N_i-k_i$ is at least $cq$, for a fixed
$c>0$ and all large $q$. This follows by bounding the error from
$qt$ in each coordinate by $\delta q\sum_\ell\|u_\ell\|_\infty$
and using the strict stock margins. Each $m_\ell$ ranges over at
least $e+1$ consecutive integers for large $q$.

If $A(k)$ vanished throughout this grid, then
$A(qb+\sum_\ell m_\ell u_\ell)$ would vanish for every
$m\in\Z^h$. Indeed it is a reduction of a finite linear combination
of rational polynomials of total degree at most $e$ taking integer
values at integers, using integer representatives for the field
coefficients. Its forward difference of order $e+1$ in each
variable vanishes. The two end coefficients in that difference
are $\pm1$, so $e+1$ consecutive zeros give zeros at every integer
in the variable, first to the right and then to the left. Applying
this argument one variable at a time extends the zero grid.
Evaluating at $m=0$ gives $A(qb)=A(0)=1$, a contradiction: all
binomials of orders at most $e<q$ are $q$-periodic modulo $p$.
If $h=0$, the sole candidate $k=qb$ already has that value.
Hence one feasible selection has quotient sum zero and lower
component whose negative occurs as a subsum of $T$. Their union
is nonempty, with the asserted linear stock margins.
\end{proof}

\begin{proposition}\label{slice-relations}
{Retain the hypotheses of Theorem~\ref{slice-main} on the kernel stock,
the class sizes, and the fixed matrix $V_J$, including the existence of
$0<t_i<\lambda_i$ with $V_Jt\in V_J\Z^J$. The actual quotient columns
on $J$ may be the columns of a moving matrix $Y_q$ over $C_q$,
provided $Y_q u=0$ for every $u\in\ker_{\Z} V_J$. The conclusions of
Theorem~\ref{slice-main} still hold.}
\end{proposition}
\begin{proof}
The quotient-zero counts used in the proof {of Theorem~\ref{slice-main}} are all of the form
$qb+\sum m_\ell u_\ell$. Their image under $Y_q$ is zero because
the group has exponent $q$ and each $u_\ell$ is an integer relation.
Every other step in the proof concerns only lower coordinates.
\end{proof}

\subsection{Lifting through a bounded cofactor}\label{sec:lift}
\begin{theorem}\label{lift}
Let $m,r$ be fixed, let {$q$ run through arbitrary positive integers
tending to infinity}, and put $n=mq$.
Suppose a zero-sumfree family of length $rn+o(n)$ in $C_n^r$ has,
after reduction modulo $q$, a basis with at least $n-o(n)$
positions in each ordinary row. Along a subsequence it has an
ordinary flag over $C_n^r$ with $n-o(n)$ positions in each row.
\end{theorem}

\begin{proof}[{Proof of Theorem~\ref{lift}}]
\textit{Preparing coordinates.}\quad
In a cyclic group $C_L$ write $H_D=\{x:Dx=0\}$; the ambient modulus
will be clear. Let $a$ be the product of the full prime powers in
$n$ at primes not dividing $q$, and put $L=n/a$. Then $a\mid m$,
$q\mid L$, $\gcd(a,L)=1$, and $L$ and $q$ have the same prime
divisors. Pass to a subsequence with $a$ constant and set $b=L/q$.
Thus $b=m/a$ is fixed and $C_n^r=C_L^r\oplus A$ where $A=C_a^r$.
The given basis modulo $q$ lifts to a basis of $C_L^r$: lift the
entries arbitrarily; its determinant is prime to $q$, hence to $L$.
The row sizes are $n+o(n)$, since they are disjoint and their total
is $rn+o(n)$.

We henceforth take subsequences as needed. For a list of scalars
in a row, call it bounded in this paragraph if, apart from $o(n)$
positions, it lies in a single coset of $H_D$ for some fixed $D$.
The coset can vary with $n$. Pass to a diagonal subsequence on
which the largest coset proportion for each positive integer $D$
has a limit. If the list is not bounded, then for each fixed $D$
at least $\delta_D L$ terms are outside a largest coset, for some
$\delta_D>0$ and all large $L$. Call this second possibility spread.
Subtracting a common scalar and then adding bounded-order errors
preserves spread: if $E$ kills all errors, concentration of the new
values in a coset of $H_D$ would put the old values in a coset of
$H_{DE}$. Deleting $o(n)$ terms also preserves spread.

Process the rows in their original order. When processing row $j$,
earlier basis vectors have been replaced by actual primary anchors
in their rows. Their reductions modulo $q$ are triangular with
diagonal one in the original basis, so the primary coordinates
in the current basis satisfy $x_j=1\bmod q$ and $x_i=0\bmod q$
for $i>j$. Thus $x_j-1$ and all $x_i$ for $i>j$ belong to $H_b$.
Classify the lower coordinates as bounded or spread. Delete the
$o(n)$ exceptions in the bounded coordinates. The bounded lower
coordinates, together with the diagonal and higher coordinates,
now have only boundedly many possible tuples, because a coset of
$H_D$ has at most $D$ elements. Choose an actual primary anchor $y$
for which this entire tuple and the cofactor mark $c_j$ are shared
by at least a fixed positive fraction of the row. This is possible
by the pigeonhole principle, also when all lower coordinates are
spread.

Replace basis column $j$ by $y$. Its $j$th entry is a unit modulo
$L$ since it is one modulo $q$ and these moduli have the same prime
divisors. In the new coordinates $x_j$ becomes
$t=x_j/y_j\in1+H_b$ and $x_i$ becomes $x_i-t y_i$ for $i\ne j$.
Every bounded lower coordinate is consequently in a bounded-order
subgroup: $x_i-y_i$ was killed by its $D$, and $(t-1)y_i$ by $b$.
Every spread coordinate remains spread. All higher coordinates are
still in $H_b$. On the positive fraction sharing the anchor tuple
and mark, the diagonal is exactly one, all bounded lower coordinates
and higher coordinates are zero, and the mark is $c_j$.
Later basis changes perturb the coordinates of this row by elements
of $H_b$ only: the coefficient at each later pivot lies in $H_b$,
division by its pivot unit keeps it in $H_b$, and multiplication
by any entry preserves $H_b$. The positive fraction just specified
has zero later coefficients and is unaffected.

After the $r$ stages take a common multiple of all torsion bounds
and pass to a subsequence on which its gcd with $L$ is a fixed
integer $B$. Then $B\mid L$, $\gcd(a,B)=1$, and all these bounded
errors are in $H_B$. For each position in row $j$ its primary
coordinates in the resulting basis have the form
$$
                 e_j+w+u .
$$
Here $w\in H_B^r$ collects all bounded deviations, including the
diagonal error and all higher components, and $u$ is supported
on the spread coordinates $i<j$; we set $w_i=0$ on these.
The joint mark is $(c,w)\in A\oplus H_B^r$. Use the isomorphism
$C_B\to H_B$ given by $t\mapsto (L/B)t\bmod L$. Pass further so the
frequency of each mark has a limit and discard marks of zero
limiting proportion. Each remaining mark has $\gg L$ positions.
The anchor mark $(c_j,0)$ survives, by its positive fraction above.
The subgroup generated by the remaining differences
$(c-c_j,w)$ is a direct product $D_j\oplus E_j$ of its two
projections: multiplication by an integer equal to one modulo $a$
and zero modulo $B$ isolates the first projection, and subtraction
isolates the second. This includes $a=1$ and $B=1$.

\medskip\noindent\textit{Prescribing a selection in one row.}\quad
Put $M=aB$, $F=\operatorname{lcm}(1,\ldots,r)$.
Suppose an integer $T$ satisfies $\eta L\le T\le N-\eta L$
in a row of size $N=n+o(n)$, where $\eta>0$ is fixed.
Given arbitrary $d\in D_j$, $e\in E_j$ and offsets on every
spread coordinate, we can select $T$ positions with joint mark
$(Tc_j+d,e)$ and leave each updated spread-coordinate residual
at distance at least $\eta' L$ from $H_B$ for a constant $\eta'>0$.
The distance here is the circular distance in $C_L$ to the nearest
element of $H_B$. The updates in this assertion use the $u$ entries.

For details, in each spread coordinate $i$ the number outside the
largest coset of $H_{MBF}$ is at least $\delta L$ for fixed
$\delta>0$; take the minimum over the finitely many coordinates and
positive mark proportions, decreasing $\delta$ if necessary.
Reserve $\lfloor\kappa L\rfloor$ disjoint choices for each such
coordinate with disjoint sites across coordinates. Choose $\kappa$
so that $2Mr\kappa<\min(\eta,\delta)/4$ and
$M\kappa<\delta/(16|A\oplus H_B^r|^2)$. After earlier reservations
at least $3\delta L/4$ positions are still outside a largest coset.
Greedily take $\lfloor\delta L/4\rfloor$ disjoint pairs whose
differences are outside $H_{MBF}$. At every step fewer than
$\delta L/2$ additional sites have been removed, so two cosets
remain. By pigeonhole retain a common ordered pair of joint marks
on at least $\lfloor\delta L/4\rfloor/|A\oplus H_B^r|^2$ pairs.
Among any $M+1$ such differences some $M$ have sum outside
$H_{BF}$. Otherwise subtracting the sums omitting different terms
makes all $M+1$ differences equal modulo $H_{BF}$; then their
common class is killed by $M$, so every original difference belongs
to $H_{MBF}$, a contradiction. Group $M$ pairs with the required
property successively; the preceding bounds give enough choices.
Only those sites are reserved. The options use all first members
or all second members. They have size $M$, equal zero joint mark,
and scalar difference $\Delta\notin H_{BF}$. Its image in
$C_L/H_B$ has order exceeding $r$, since its order cannot divide
$F$.

Select one option in each choice. Write $G$ for the resulting
number of selected sites, a multiple of $M$, with $2G$ sites
reserved. For the other $T-G$ sites use the mark classes among the
$N-2G$ unreserved positions. The ideal counts proportional to class
sizes are a fixed positive multiple of $L$ away from both endpoints
in every class, by the bounds above. In the group of mark
differences express $(d,e)$ as an integral combination of the
differences from $(c_j,0)$. Round every nonbase class to its
prescribed coefficient modulo $M$, within $M$ of its ideal count.
Choose the base count to make the sum $T-G$. Its error is at most
$M|A\oplus H_B^r|$ and is absorbed by the available margin. The
joint mark sum is $(Tc_j+d,e)$ because $M$ kills the finite mark
group and $G$ is a multiple of $M$.

{Fix a spread coordinate $i$ and use the collection reserved for $i$.
Let $Y$ be the set of its attainable scalar sums in coordinate $i$
modulo $H_B$.} The cyclic quotient has order $L/B$.
If the stabilizer of $Y$ has size $\le r$, no choice difference
can stabilize any intermediate sumset, because its order exceeds
$r$ and a period persists under addition. Thus $|Y|\ge
\lfloor\kappa L\rfloor+1$. If its stabilizer has size $h>r$, it
is a union of cosets with spacing $L/(Bh)$ in the quotient circle.
Use the uniform distribution on $Y$ by choosing a representative
arrangement for each sum. Across collections make these choices
independently. Conditional on all other collections, the probability
of distance less than $\eta' L$ from $H_B$ is at most
$$
 \max\left\{\frac{2\eta'L+1}{\lfloor\kappa L\rfloor+1},
                 \ 2B\eta'+\frac1{r+1}\right\}.
$$
This follows by counting the points in a circular interval of
length $2\eta'L$ in $C_{L/B}$, first by all integer positions and
then by equally spaced points of each stabilizer coset. Choose
$\eta'$ sufficiently small that the bound is less than $1/r$.
The union bound over at most $r$ spread coordinates gives all
desired inequalities simultaneously. With no spread coordinate the
mark selection alone suffices.

\medskip\noindent\textit{The obstruction from a directed cycle.}\quad
Draw a spread edge $j\to i$ for each spread coordinate of row $j$.
All spread edges point downwards. Draw a torsion edge $j\to i$ if
some vector of $E_j$ has nonzero $i$th coordinate. Such an edge
can go in either direction or be a loop.

We describe how to select actual indexed positions. Prescribe
count residues $k_j\bmod a$, cofactor deviations $d_j\in D_j$,
and torsion deviations $e_j'\in E_j$. Put the sum of all prescribed
torsion deviations into the initial residual $z=\sum_j e_j'$
in $C_L^r$. The procedure then updates only with the diagonal
counts and the spread entries $u$. This accounts for torsion
exactly once if each prescribed deviation is realized. Process
rows in descending original order. In row $j$ take the unique
$T\in[0,n)$ with
$$
                    T=-z_j\bmod L,\qquad T=k_j\bmod a .
$$
If some previously {active} row has a spread edge to $j$, the
current $z_j$ has distance at least a positive multiple of $L$
from $H_B$. This holds by induction: whenever a spread update
changes that coordinate, use the preceding one-row construction
with the current offset, and between such updates it remains fixed.
Otherwise $z_j$ is its initial value in $H_B$. Each nonzero
element of $H_B$ has distance at least $L/B$ from zero. Therefore
unless $z_j=0$ and $k_j=0$, both $T$ and $n-T$ are at least a
positive multiple of $L$. In particular for $z_j=0$, $k_j\ne0$,
$T$ is a positive multiple of $L$ below $n$. The row size $n+o(n)$
provides the needed stock margin in each case. Select $T$ sites
using the preceding construction with deviations $d_j,e_j'$,
cancel the diagonal and update all lower spread coordinates.
Call such a row active. When $T=0$ skip the row; this step is valid
provided both prescribed deviations there are zero. Only $r$
updates occur, so successive positive margins remain positive.

If the $c_j$ fail to give a basis for $A$, the map
$C_a^r\to A$, $k\mapsto\sum k_jc_j$, has a nonzero kernel element
by equality of the finite orders. Taking this $k$ and zero
deviations throughout gives an indexed zero: every row with
$k_j\ne0$ is active and the whole cofactor sum is zero. Hence on
a zero-sumfree family we use the $c_j$ as cofactor basis (and
ignore cofactor coordinates if $a=1$). Add a cofactor edge
$j\to i$ if any $d\in D_j$ has nonzero $i$th coordinate.

Suppose the union graph has a directed cycle. Choose a shortest
one, assigning each of its edges one of its available types.
Since spread edges alone point downwards, at least one assigned
edge is of torsion or cofactor type. At each torsion source choose
$e_j'$ nonzero at its cycle successor, and at each cofactor source
choose $d_j$ nonzero there. Set all other deviations to zero.
Take $k$ equal to the negative coordinate vector of $\sum_jd_j$
in the cofactor basis. At every cofactor-edge destination on the
cycle $k_j\ne0$, and at every torsion-edge destination the initial
$z_j\ne0$. For each assertion the contributing predecessor is
unique, since any edge from another cycle vertex to that destination
would yield a shorter directed cycle. For a loop there is a single
source. Initial residuals are defined using the prescribed shifts
even for sources later in the processing order.

Each such destination is active in the selection procedure: a
spread update, if it occurred, left its residual far from $H_B$;
otherwise its specified nonzero initial residual or nonzero count
residue suffices. Follow the cycle from these destinations along
successive spread edges. Every such edge goes downward in processing
order, and a selected source makes its target active. Thus all
sources of nonzero prescribed deviations are active, so no skipped
row was required to realize a nonzero deviation. Formally this is
an induction in descending processing order: a source is either
itself a torsion or cofactor destination, or has an earlier-processed
active predecessor on such a spread segment. All deviations are
therefore realized and at least one row contributes a position.
The final primary coordinates vanish by their diagonal corrections;
the cofactor sum is $\sum k_jc_j+\sum d_j=0$. This is a nonempty
zero in the original group.

The union graph must consequently be acyclic on a zero-sumfree
family. Reorder the anchors with every edge pointing to an earlier
vertex. Pair each primary anchor with its cofactor base mark $c_j$;
these pairs form a basis of $C_n^r$ by the coprime decomposition at
$L,a$. Every retained row vector has its anchor with coefficient
one and only earlier deviations: the spread, torsion and cofactor
edges account for all possible deviations, and acyclicity excludes
diagonal loops. Thus each full ordinary row contains $n-o(n)$
positions. {This proves Theorem~\ref{lift} on the chosen subsequence.}
\end{proof}

\subsection{Ordinary stability at unrestricted exponents}

\begin{theorem}\label{stability}
Fix $r\ge1$. Every family of zero-sumfree sequences in $C_n^r$
with $n\to\infty$ and $|S_n|/n\to r$ has a subsequence admitting
bases with $n+o(n)$ indexed positions in each ordinary row.
\end{theorem}

\begin{proof}[{Proof of Theorem~\ref{stability}}]
\textit{Passing through a kernel at a growing prime power.}\quad
Use Theorem~\ref{chart-main}. In a group
$G=C_n^d\oplus E$, $|E|=O(1)$, $nE=0$, it gives a map
$F:G\to C_n^k$ of bounded image index, with kernel
$C_n^\ell\oplus E'$ for $\ell=d-k$ and $|E'|=O(1)$, containing
$\ell n+o(n)$ positions. Off the kernel there are finitely many
fixed positive-mass output values in $\Z^k$ capturing $k n+o(n)$
terms, partitioned into strictly balanced lists in complementary
subspaces.

Choose free coordinates $C_n^d=U\oplus L$, where $U=C_n^\ell$ is
the direct summand of $\ker F$ obtained by Smith reduction on the
free part {in the proof of Theorem~\ref{chart-main}} and $L=C_n^k$. Thus $G=U\oplus(L\oplus E)$ and
$E'=(L\oplus E)\cap\ker F$. {Indeed $\ker F=U\oplus E'$ and
$|\ker F|=n^{d-k}|E|[C_n^k:\operatorname{im} F]$, so this particular
complement has bounded order.} Write $B$ for a fixed common multiple
of the exponents of $E,E'$, passing to a subsequence to take them
constant. We may use their least common multiple, so $B\mid n$.
Within each fixed output $v$, including $v=0$, split positions
according to their complete $L\oplus E$ component $w$. There
are at most $|E'|$ such components since $U\le\ker F$. Discard
subtypes of zero limiting mass after stabilizing all frequencies.
The retained number is finite. Every integral relation among the
labels $v$ gives, on the corresponding selected exact labels $w$,
a sum in $E'$.

Let $q=p^a\to\infty$ be any prime power dividing $n/B$ and put
$m=n/q$. {For $\ell\ge1$, let $\mathsf s(C_m^\ell)$ denote the least length
forcing a zero-sum subsequence of exactly $m$ indexed terms. By
Alon and Dubiner~\cite{alon-dubiner}, {there is a constant $c_\ell$ depending only on $\ell$ such that}
$\mathsf s(C_m^\ell)\le c_\ell m$. Apply the bound in each subtype
to the $U$-components reduced modulo $m$.} Repeat
greedily, leaving $O(m)$ unused terms per subtype. For $\ell=0$
or $m=1$ just group directly. Each block belongs to
$m C_n^d$ in its free coordinates and has zero $E$ component
because $B\mid m$. Divide its sum by $m$ in the canonical
identification $m C_n^d\simeq C_q^d$. {The map $mx\mapsto x\bmod q$
is well defined because $mx=mx'$ modulo $n=mq$ exactly when
$x=x'$ modulo $q$. It works even when $m$ and $q$ have a common prime
factor.} In $L$ its new value is
$w_L\bmod q$. A relation of the preceding paragraph has
$L$ component killed by $B$ in $C_n^k$ and hence divisible
coordinatewise by $n/B$, a multiple of $q$. Thus all integral
relations among the fixed $v$ are respected by the divided upper
values, which are equal even across subtypes of the same $v$.
All zero-label blocks lie in $C_q^\ell$ and there are
$\ell q+o(q)$ of them. For output $v\ne0$ of mass $\lambda_v$
there are $\lambda_vq+o(q)$ blocks. These assertions follow by
division by $m$ with an additive error $O(1)$ per subtype. The
block sequence is zero-sumfree since each selection expands to
disjoint nonempty blocks.

{Apply Proposition~\ref{slice-relations}, the moving-column form of Theorem~\ref{slice-main}, to the block sequence.} Its kernel stock lies in the required
direct summand $C_q^\ell$, has length $\ell(q-1)-e$ where
$0\le e=o(q)$ {by Lemma~\ref{slice-residual} and the count of kernel blocks}, and
each upper class has the required positive stock. {The relations verified above are exactly those required of the upper columns by Proposition~\ref{slice-relations}.} Thus a strict
selection on any of our fixed output sublists gives a nonempty
actual zero, a contradiction. There always is a growing prime
power divisor of $n/B$: otherwise the prime divisors and all their
valuations would be bounded. This proves absence of a strict
selection at every stage, with no bound on the cofactor $m$.

{A constituent of dimension at least two satisfies
Theorem~\ref{geom-main} in the lattice generated by its columns.
That theorem would give a strict selection, which the preceding
argument excludes. Every constituent at this stage therefore has
dimension one.}

\medskip\noindent\textit{Rank-one constituents and primes of bounded size.}\quad {On each such line all occupied columns have a single oriented value.} Indeed opposite signs give a positive real
relation of arbitrarily small coefficients. For columns of the
same sign, choose the generator of their integer span for which all
the columns have positive coefficients, and write them in that basis. If any exceeds one, their weighted sum exceeds
one since the masses sum to one. Multiplying all capacities by
the reciprocal of that weighted sum gives strictly interior
coefficients attaining that generator, an integer combination of the columns. Thus they
are all the same value. Denote the resulting $k$ independent
integer columns by $v_1,\ldots,v_k$; each has mass one. Let $w_j$
be one actual upper label of positive density above $v_j$, and
write $a_j=(w_j)_L$. Every other subtype above it differs in $L$
by an element killed by $B$. Zero-label subtypes have the same
bounded-torsion property about zero.

We show that $a_1,\ldots,a_k$ form a basis of $C_{n/M}^k$ after
reduction, for a fixed divisor $M$ of the exponents on a subsequence,
with $B\mid M$. Write $A$ for their column matrix and $F_L$ for
the matrix of $F$ on $L$. Since $B E=0$,
$$
              B F_L A\equiv B [v_1,\ldots,v_k]\pmod n .
$$
Consequently a prime dividing both $n$ and $\det A$ belongs to the
fixed finite set of primes dividing $B\det[v_1,\ldots,v_k]$.
For each of these primes pass to a subsequence where its valuation
in $n$ is constant or tends to infinity, and where divisibility
of $\det A$ by it is constant. If $p$ is such a prime of unbounded
valuation, take $q$ equal to the full $p$-part of $n/B$ in the
preceding construction. Suppose $A$ is singular modulo $p$.
Choose $b\ne0$ modulo $p$ with $A b=0\bmod p$, with
$0\le b_j<p$; $p$ is fixed. For each $b_j>0$, prescribe exactly
$c_j=b_jq/p$ divided upper blocks; put $c_j=0$ at other classes.
These counts are admissible after the $O(e+1)=o(q)$ trimming {of Lemma~\ref{slice-trim}}, since every positive $b_j/p$ is strictly
between zero and one and every class has $q+o(q)$ blocks.
Their upper sum is $(q/p)A b=0\bmod q$.

{Because the exceptional prime $p$ is fixed, $e=o(q)$ implies
$e<q/p$ for all sufficiently large $q$. Every prescribed count is a
multiple of the prime power $q/p$. By the {\hyperref[lem:prescribed-counts]{lemma on correction at prescribed counts}}, its count function $H$ satisfies $H(c)=H(0)=1$, and some
selection with those counts has a lower component canceled by the
kernel stock. If $\ell=0$ that correction is empty. The quotient sum is
zero and at least one $c_j$ is positive, giving a nonempty zero. {This contradicts the zero-sumfreeness of the block sequence. Hence $A$ is nonsingular modulo $p$.}}

Take $M$ to be the least common multiple of $B$ and the full,
constant prime powers in $n$ from the exceptional primes of
bounded valuation. The determinant of $A$ is a unit modulo
$n/M$ by the preceding argument. Also every $w$ above $v_j$
has $L$ component $a_j\bmod(n/M)$, while the zero group has
$L$ component zero there. This proves a genuine basis and $k$
independent cosets above the lower free coordinates after division
of the modulus by a bounded factor.

\medskip\noindent\textit{{Iteration.}}\quad
Apply the same construction to the zero-label positions as a
sequence in $U\oplus E'$. They have length $\ell n+o(n)$ and
are zero-sumfree. All subsequent choices of free coordinates are
made \emph{inside $U$}. The earlier prescription of the upper
$L$ components is therefore compatible with each recursion.
{Each stage passes from free rank $d$ to $d-k$ with $k\ge1$.
Iteration stops when no free coordinate remains, after at most $r$
stages.} Take the least common multiple $M_*$ of the bounded
divisors chosen at those stages, after all subsequences have been
taken. Then $M_*\mid n$ is fixed. Each basis at a stage remains
a basis modulo $n/M_*$, a divisor of its earlier modulus.
Choose the bases in increasing order, starting at the deepest
kernel and adjoining each stage's $k$ independent upper basis
vectors. For an occupied term at such a stage all further upper
coordinates vanish modulo $n/M_*$, and its current $k$-tuple is
one of the indicated basis vectors. Lower coordinates are free.
Thus its output lies in the associated ordinary flag row, in any
order of these $k$ vectors. Each of the $r$ rows contains at
least $n-o(n)$ actual positions after reduction modulo $n/M_*$.
Their total is $rn+o(n)$, so each also has at most $n+o(n)$.

Apply Theorem~\ref{lift} with $m=M_*$ fixed and $q=n/M_*\to\infty$.
It gives the desired flags at the original modulus. The upper bounds
follow by subtracting the lower bounds in the other disjoint rows
from $rn+o(n)$. This proves Theorem~\ref{stability}.
\end{proof}

\section{Preparation of the flag rows and matching of their counts}\label{sec:tools}

In this section we give the estimates for selecting indexed positions
inside a row. We also partition the rows so that each lower coordinate
is either constant on a class or can be varied by exchanging disjoint
pairs of positions. In Subsection~\ref{sec:counts} we use this
partition to match the counts of classes in an auxiliary sequence
to counts in the original sequence.

\subsection[Cyclic estimates]{{Cyclic estimates}}

From now on suppose a critical zero-sumfree family has ordinary flags
with designated populations $n-1-d_i$ and exterior size
$k=1+\sum d_i\to\infty$, $k=o(n)$. We need two estimates for
coordinating symbolic selections on actual positions.

Write $\norm{x}_n$ for circular distance to zero in $C_n$.
\begin{lemma}\label{cyc-prog}
Let $k\to\infty$, $k=o(n)$ and $v_n\in C_n$. If
$\norm{h v_n}_n\to\infty$ for every fixed integer $h\ge1$, then,
on any interval of $\Theta(k)$ consecutive integers $j$, only $o(k)$
values of $a+jv_n$ fall in a given interval of at most $k$ consecutive
residues, uniformly in the translate $a$ and in the initial index.
\end{lemma}
\begin{proof}
If a fixed positive proportion hit on a subsequence, two successive
hits would have bounded index gap, since the sum of their gaps is
$O(k)$. Pass to a subsequence where one such gap is a fixed $h$.
Then $b=\norm{hv_n}_n\le k$ and $b\to\infty$. Split indices into
their $h$ residue classes. In each, lift to a real progression whose
step is $b$ or $-b$. The total traversed distance is $O(kb)$, so
it can meet $O(kb/n+1)$ translates by $n$ of the tested interval.
Each translate contains $O(k/b+1)$ points. Their product divided
by $k$ is $O(k/n+b/n+1/b+1/k)=o(1)$. Summing over the fixed $h$
contradicts the asserted positive proportion.
\end{proof}

\begin{lemma}\label{cyc-sw}
For $l=\operatorname{lcm}(1,\ldots,r)$ and $k=o(n)$, if the optional
sums of $v_1,\ldots,v_m\in C_n$ do not contain $r+1$ residues at
pairwise distance greater than $k$, then
$\sum_j\min(k,\norm{lv_j}_n)\le(20rl+r(r-1))k$ for all large $n$.
\end{lemma}
\begin{proof}
A maximal separated set covers all optional sums by at most $r$ balls
of radius $k$. Assign each increment to a ball containing it. In a
ball of center $a$ containing at least $r$ assigned increments, order
any $r$ of them. Two of their $r+1$ prefix sums lie in a common ball.
Their difference is within $2k$ of zero and is a sum of $q$ of the
assigned increments, $1\le q\le r$. Thus $\norm{qa}_n\le(q+2)k$
and $\norm{la}_n\le3lk$. Every assigned increment in that ball has
$\norm{lv_j}_n\le4lk$. At most $r(r-1)$ increments are in the
other balls. Ignore them for now. The powered optional sums are
covered by at most $r$ balls of radius $lk$. Follow the positive
centered representatives of the remaining powered increments by adding
them successively. Every prefix is covered, and each step has length
at most $4lk$. The entire continuous path lies in the same balls
enlarged to radius $5lk$, of total length at most $10rlk<n$.
It cannot circle the group. Its total length is at most $10rlk$.
Repeating on the negative representatives and restoring the ignored
increments gives the bound.
\end{proof}

\subsection{Preparing the triangular stock classes}\label{sec:prep}
Write $C_n=\mathbb Z/n\mathbb Z$ and $\norm{a}_n$ for circular
distance from $a\in C_n$ to zero. Let $r$ be fixed, $n\to\infty$,
$k\to\infty$, and $k=o(n)$. In each row of an ordinary flag of
$C_n^r$ retain $n-1-d_j$ indexed positions with $0\le d_j=O(k)$.
The ordinary row of a basis vector is its translate of the span
of all preceding basis vectors.

\begin{lemma}\label{lem:prepare-stock}
After passage to a subsequence, a triangular change of basis and
removal of $O(k)$ positions, each row is partitioned into a fixed
finite number of indexed classes. Every class has size $\omega(k)$.
For each class and each lower coordinate $i$, one of these two
assertions holds:
\begin{enumerate}
\item The coordinate is constant within the class at each $n$;
its value may depend on $n$.
\item {For all fixed integers $h\ge1$ and $H\ge0$}, the maximum number
of disjoint pairs $x,y$ in the class with
$\norm{h(x_i-y_i)}_n>H$, divided by $k$, tends to infinity.
\end{enumerate}
Each row has a distinguished class with at least $c n$ positions
for some fixed $c>0$ whose constant lower coordinates all vanish.
One of its actual positions is the new basis vector of the row.
\end{lemma}
\begin{proof}
First prepare the classes in a single row in a given coordinate
system. At each stage pass to a subsequence on which the size of
the current class divided by $k$ tends to a finite limit or to
infinity. In the former case remove this class. In the latter,
if every coordinate outside those already made constant satisfies
the second assertion, retain the class. If the assertion fails
for some coordinate $i$, some fixed $h,H$ and a subsequence have
maximum matching size at most $Dk$ for a fixed $D$. Discard the
positions of that matching. For any two remaining positions the
powered difference in coordinate $i$ has norm at most $H$, for
otherwise the matching could be enlarged. Relative to one
remaining value $a$, every value $v$ then satisfies
$h(v-a)\in[-H,H]\pmod n$. Each of these at most $2H+1$ images
has at most $\gcd(h,n)\le h$ preimages. Thus there are at most
$h(2H+1)$ values. Pass to a fixed number of nonempty values,
order them by multiplicity, and split by the actual values.
The coordinate just tested is now constant on every child and
remains so under further restrictions.

Run the same argument on every child, retesting all coordinates
that have not been fixed. A property of dispersion before
splitting is not needed after splitting. Every split fixes one
more coordinate and each split has finite branching. Induction on
the number of coordinates therefore gives a finite tree and
removes only $O(k)$ positions in total. Later subsequences
preserve all limits already obtained on the other branches.
Constant coordinates which happen incidentally are included in
the first assertion.

Now treat the rows in increasing order. At row $j$, express its
lower part using the already chosen basis vectors and apply the
preceding partition. Its original size is $n-O(k)$, and after
discarding $O(k)$ positions its finitely many classes still
account for $n-o(n)$ positions. One of them has size at least
$c_j n$ along a subsequence for some $c_j>0$. Choose a position
of that class as the new representative $e'_j$. Its coefficient
at the old $j$th basis vector is one, so this operation is an
invertible unit triangular change of basis. On row $j$ it
subtracts a common vector from all lower components. Pair
differences are unchanged, and all constant coordinates of the
distinguished class become zero. The individual ordinary rows
of the flag and the coordinates already fixed on earlier rows
are preserved. Complete all $r$ steps and take $c=\min_j c_j$.
\end{proof}

If the initial number of exterior positions was $k$, reclassifying
the removed positions makes it $k'$, where $k\le k'=O(k)$. The
displayed dispersion assertions and all the class size lower
bounds persist with $k'$ as scale. They also provide disjoint
collections of $Ak'$ pairs in any fixed finite list of the
dispersed coordinates, for any fixed $A>0$, at any prescribed
fixed thresholds and multipliers. To see this, select collections
in succession; each previously occupied position can invalidate
at most one edge in a matching for the next request.
This preparation gives the class decomposition used in a
deficit-scale construction. Transferring a zero requires also
choosing compatible numbers of terms from these classes when
different rows contribute to the same lower coordinate.

\subsection[Exterior patterns and exchanges with prescribed counts]{{Exterior patterns and exchanges with prescribed counts}}\label{sec:exterior-exchanges}

{
\medskip\noindent\textit{The exterior distribution.}\quad
Use the basis and row deficits obtained from Lemma~\ref{lem:prepare-stock}.
From now on write $k$ for the enlarged exterior size.
Pass to a subsequence with $d_i/k\to\delta_i$, so that
$\sum_i\delta_i=1$. In each coordinate order the distinct scalars
of the exterior by decreasing frequency, breaking ties in any fixed
way. Pass to subsequences so every fixed frequency divided by $k$
has a limit $p_a$. Here and below the labels $a$ and the limits refer
to the given coordinate. Their sum is at most one. At stage $n$ keep
the names $a\le J_n$ having $p_a>0$, with a cutoff $J_n\to\infty$
slow enough that the total {absolute} error in the {normalized} frequencies up to $J_n$,
summed in each coordinate, tends to zero. Such a choice follows by
first imposing any fixed cutoff and then increasing it diagonally.
Use a star for every other coordinate value. A star's largest single
scalar multiplicity is $o(k)$: a positive limiting frequency at
any fixed place is named eventually, and positions past any fixed
cutoff $A$ have frequency at most $1/(A+1)$ by their ordering.

Each exterior position now has a pattern of $r$ entries. Pass again
so each fixed pattern $t$ has limiting mass $\mu_t$. These masses
sum to one. For each coordinate the total named mass with labels
above $A$ is at most $\sum_{a>A}p_a+o(1)$ by the choice of $J_n$;
a union bound over the coordinates gives tightness on this countable
pattern space. For each fixed name pass diagonally so its equality
tests with zero and one eventually stabilize.

For any finite set $F$ of positive-mass patterns use the following
description of evaluation on its named entries.

For coordinate $i$ let $F_i$ consist of the patterns in $F$ named
there. Let $f_{i,n}:\Z^{F_i}\to C_n$ be evaluation at their true
coordinate values. Pass to a subsequence so, for every integer
direction $u$ and positive $h$, the centered representatives of
$hf_{i,n}(u)$ are eventually constant or have absolute value tending
to infinity. Define $T_i$ to consist of those $u$ with a bounded
multiple. This is a subgroup: common multiples handle addition and
negation. Choose an integral basis of $T_i$ and a common positive
$L$ for all its basis multipliers and all $i$. On each $T_i$ there
is an integral homomorphism $A_i$ such that eventually, on the
whole subgroup simultaneously,
$$
                 Lf_{i,n}(u)=A_i(u)\pmod n .
$$
Indeed define $A_i$ by the eventual centered integer values on its
basis, then extend linearly. There is also a phase homomorphism
$\phi_{i,n}:T_i\to\Z/L\Z$ given by
$(L\widetilde f_{i,n}(u)-A_i(u))/n\pmod L$, using any integral
representative of the evaluation. Changing the representative adds
a multiple of $L$ to the quotient. Pass to a further subsequence
so its finitely many basis values stabilize. Denote the resulting
phase by $\phi_i$. The assertions now hold for every $u\in T_i$
at every sufficiently late $n$, even if $u$ is chosen after $n$ is fixed.
\par}

{
\medskip\noindent\textit{Switches with a prescribed count.}\quad
In a class $A$ of size $\omega(k)$ call a lower coordinate $i$
dispersed if for each fixed $t,H\ge1$ it has $\omega(k)$ disjoint
pairs with $\|t(a_i-b_i)\|_n>H$. {Lemma~\ref{lem:prepare-stock}} gives
classes in which every lower coordinate is either constant or
dispersed. Let $W_A$ be the latter set.

\begin{lemma}\label{pos-banks}
Fix $H>0$ and $\epsilon>0$. In finitely many disjoint classes,
prescribe {an integer} selected count $z_A$ {with $0\le z_A\le |A|$ in every class}. For each class to be used in an
exchange assume $\epsilon k\le z_A\le |A|-\epsilon k$.
For each of a list of at most $r$ distinct target coordinates
choose one such class dispersed in that target. Given an arbitrary
translate and a forbidden interval of length at most $Hk$ in each
target, {for all sufficiently large $n$} there is a simultaneous selection of exactly $z_A$ positions
per class avoiding all these intervals in the sums. The contributions
in every coordinate constant on each of the assigned classes are
independent of the exchanges.
\end{lemma}
\begin{proof}
Write $l=\operatorname{lcm}(1,\ldots,r)$. At each assigned coordinate
reserve $\lfloor\gamma k\rfloor$ disjoint pairs in its class, with
$2r\gamma<\epsilon$. The powered differences in its target may all
be taken to have norm at least $H_n$ for some $H_n\to\infty$:
diagonalize in the fixed threshold, then discard all pairs meeting
a previously reserved site. Each matching loses at most as many
edges as the number of sites already reserved, which is $O(k)$.
For each reserved pair select exactly one endpoint. In class $A$
with $b$ such pairs there are $|A|-2b$ other sites, from which
we must take $z_A-b$. Both this number and the number left unused
are nonnegative because $b\le\epsilon k$. Take that many fillers
and hold them fixed.

Call the pairs reserved for one target its bank. Write $L=\lceil Hk\rceil+2$. Within a bank the sum of
$\min(L,\|l(a_i-b_i)\|_n)$ over its pairs is at least
$\lfloor\gamma k\rfloor\min(L,H_n)$, which divided by $L$
tends to infinity. The cyclic switch estimate in Lemma~\ref{cyc-sw} gives $r+1$ configurations of the whole bank with
target sums at pairwise distance greater than $L$. Choose uniformly
and independently from these configurations for each bank. After
fixing the other choices the probability of a hit in its target
is at most $1/(r+1)$. The sum of the at most $r$ hit probabilities
is less than one. Some configuration avoids all intervals. Each
exchange occurs within its original class and preserves all of
that class's constant entries.
\end{proof}
\par}

{
\medskip\noindent\phantomsection\label{construction:exterior-exchanges}\textit{Exchanges within exterior patterns.}\quad
The conclusion of Lemma~\ref{pos-banks} also holds when some assigned
classes are exterior patterns of fixed positive mass and their target
coordinates are starred. For each such pattern $P$ of mass $\mu_P>0$, let $P_n$ be its set of
indexed positions, so $|P_n|=(\mu_P+o(1))k$. These sets and the stock
classes are pairwise disjoint. Assume the same count margins
$\epsilon k\le z_P\le |P_n|-\epsilon k$, with $\epsilon>0$ fixed as
$n\to\infty$.

Fix a positive integer $u$ and a threshold $H_0\ge0$. For any scalar
$a$ there are at most $u(2\lfloor H_0\rfloor+1)$ values $b$
with $\|u(a-b)\|_n\le H_0$. Each such value occupies $o(k)$
positions in the starred coordinate. Choose a fixed $\gamma>0$
with $2r\gamma<\epsilon$ and $2r\gamma<\mu_P/2$ for every
exterior pattern in the finite list. Even after reserving the
sites for earlier targets and for the first pairs of the current
target, at least $(\mu_P/2+o(1))k$ sites of $P_n$ remain.
Any first endpoint among them therefore has an available second
endpoint with $\|u(a-b)\|_n>H_0$ for all large $n$. Select
$\lfloor\gamma k\rfloor$ disjoint pairs per target in succession.
As this works for each fixed $H_0$, the threshold can tend to
infinity sufficiently slowly. Take $u=\operatorname{lcm}(1,\ldots,r)$.
The count margins allow the choices of endpoints and fillers in
the proof of Lemma~\ref{pos-banks}. The same separated configurations
and conditional probability bound apply jointly to the stock classes
and the exterior patterns. Exchanges within an exterior pattern
leave its named entries unchanged.
\par}

\subsection{Joint counts in successive rows}\label{sec:counts}
We first carry out the count calculation when each stock class has one
label. The calculation specifies, for each row, which coordinates
are determined by exact equations and which can be covered by exchanges.
In Subsection~\ref{sec:realize} we use these same tests with varying entries. At a
target with {a coordinate that varies within an active class of the original sequence} we assign a bank
before applying any finite-label equation at that target. In the
finite-label case the following construction gives {zero-sumfree} auxiliary
sequences, as soon as the original critical sequences are zero-sumfree.

\medskip\noindent\textit{Compatible lattices.}\quad
{The next choice of anchors applies to the case of one label per stock class.}
Pass to a subsequence, choosing a most frequent stock label as the
basis anchor in every row. This is a unit triangular change of
basis. Reclassify all finitely many types with stock $O(k)$ as
exterior, after arranging for each ratio to $k$ to converge in
$[0,\infty]$. The new $k$ is comparable to the old one. Now every
type has stock $\omega(k)$.
{For the classes prepared by Lemma~\ref{lem:prepare-stock}, retain its
basis and distinguished classes. Their sizes already satisfy the
required bounds. In both cases measure the exterior in the chosen
basis after reclassification.}
Use $I_j$ for the coordinates indexing the nonanchor stock types
in row $j$, $s_j=|I_j|$, and $I=\coprod_j I_j$. Write $f\in\Z^I$
for the nonanchor counts, $f/n\to\beta$, and $d_i/k\to\delta_i$.
For a column of nonanchor type $a\in I_j$ write its element
$e_j+b_a$; the vector $b_a$ has nonzero coordinates only below $j$.
Use $b_{i,n}(z)$ for the $i$th scalar contribution from the lower
parts of a nonanchor count vector $z$.

After passage to a subsequence every scalar entry has either
fixed finite order or order tending to infinity: this follows by
stabilizing each test for annihilation by a positive integer.
After a further diagonalization each $u\in\Z^I$ has
$u(f-n\beta)/k$ converging in the extended real line. Define
$$
  \mathcal B=\{u\in\Z^I:u\beta\in\Q,\quad
                  u(f-n\beta)/k\text{ has a finite limit}\}.
$$
This is a saturated subgroup: addition and negation preserve both
conditions, and if a nonzero integer multiple satisfies them then so
does its divisor. For a union $P$ of label blocks put
$\mathcal B_P=\mathcal B\cap\Z^P$, taking vectors supported on $P$,
and $W_P=\{w\in\Z^P:uw=0\ (u\in\mathcal B_P)\}$.
The image of $W_{P'}\otimes\Q$ under restriction to $P\subseteq P'$
is exactly $W_P\otimes\Q$. Indeed the functionals annihilating
the image are exactly those on $P$ whose extension by zero belongs
to $\mathcal B_{P'}\otimes\Q$. Clearing denominators identifies
that set with $\mathcal B_P\otimes\Q$. Taking annihilators again
proves equality. Thus an integral $w\in W_P$ lifts after multiplication
by a fixed positive integer to an integral member of $W_{P'}$.

Write $B$ for a matrix of basis rows of $\mathcal B$, and
$\lambda=\lim B(f-n\beta)/k$. Fix rational $\rho$ with
$B\rho=B\beta$ and write $t=D\rho\in\Z^I$, $D\ge1$ fixed.
Stabilize $n\bmod D$ and every standard representative of each
$b_{a,i}$ modulo $D$. There exist fixed $C>10r+10$ and $v\in\mathbb R^I$
with $Bv=CB\beta+\lambda$, $v_a>10$, and
$C-\delta_j-\sum_{a\in I_j}v_a>10$ for all $j$. Indeed take
$C=n/k$ at a sufficiently late exponent and correct $f/k$ by $o(1)$
with a fixed right inverse of $B$; the anchor stocks give the last
margins. Every constraint in every $\mathcal B_P$ is then satisfied
by the restricted $v_P,\rho_P$.

For each of the finitely many block unions $P$ and each coordinate
$i$, stabilize the bounded-multiple tests of $b_{i,n}(w)$ for all
fixed $w\in W_P$. Either some $w$ has every fixed positive multiple
of divergent centered magnitude (a wild witness), or on $W_P$ there
is a fixed integral form $V$ with $q b_{i,n}(w)=V(w)\pmod n$.
Take a common $q$ using the finite bases and lists, and extend $V$
integrally to $\Z^P$. This is possible since $W_P$ is a direct
summand (its defining matrix has free abelian image). Choose the zero
extension if $V|_{W_P}=0$. If $V|_{W_P}\ne0$, then
\begin{equation}\label{mixdrift}
V\beta_P\notin\Q
       \quad\hbox{or}\quad |V(f_P-n\beta_P)|/k\longrightarrow\infty,
\end{equation}
since otherwise the extension by zero belongs to $\mathcal B_P$.
For $\kappa\in qW_P$ one has exactly
$b_{i,n}(\kappa)=V(\kappa)/q\pmod n$.
Also stabilize all bounded-order stock scalar phases, increasing
fixed moduli as necessary. Each such scalar is $n\theta$ for a
fixed rational $\theta\bmod1$ at every late exponent.

\medskip\noindent\textit{{Auxiliary} counts and the descending classification.}\quad
Take the exterior distribution of named scalar entries and diffuse
stars in the exterior construction of {Subsection~\ref{sec:exterior-exchanges}}. {For a finite truncation $F$, choose $\eta>0$ with $(|I|+3)\eta<1/4$. {Choose $M$ sufficiently late on the further subsequence where the named-form remainders and phases described below have stabilized} and set $K=M/C$.} Use
$\lfloor K(\mu_t-\eta)\rfloor$ copies of each positive pattern of
mass $\mu_t>\eta$ in $F$, with actual names at $M$ and stars encoded by
$2$.
Put $f'_a=\lfloor K(v_a-\eta)\rfloor$ and use total stock
$$
 N_j^*=M-1-d_j^*,\qquad d_j^*=\lceil K(\delta_j+(s_j+3)\eta)\rceil
$$
in row $j$. The anchor fills the remaining capacity. Stabilize the
named-form subgroups for integral evaluations on exterior names
and stock scalar representatives, in $C_n$ and $C_{Dn}$, as in
the subgroup argument in {Subsection~\ref{sec:exterior-exchanges}}. {Take $l$ common to all finite lists and divisible by $D,q$ and
all bounded-order denominators.} {For each named form with a bounded multiple, let $A$ and $\phi$ denote its stabilized integer remainder and phase for this choice of $l$.} {Let $m$ be the order of the ambient group of the form, so $m=n$ or $m=Dn$.} Its evaluation at each sufficiently late $n$ is
{$(A+m\phi)/l$ modulo $m$}, simultaneously on the whole subgroup,
by fixing basis remainders and phases. {Here $\phi$ denotes an integer representative of the phase modulo $l$, and the division is exact in $\mathbb Z$.} Outside it, the progression
estimate applies. This permits directions depending on a fixed $M$.

Fix a putative {auxiliary} zero with exterior counts $x_t$, nonanchor
counts $c^*$ and row totals $T_i^*$. The exterior is nonempty since
the highest row of a stock-only nonempty selection would have its
pivot count strictly between zero and $M$. Seek $hx_t$ exterior
sites, with
\begin{equation}\label{mixmult}
(1-\epsilon)k/K\le h\le k/K,\qquad h\equiv1\pmod l .
\end{equation}
The small fixed $\epsilon>0$ can depend on $M,F,\eta$. It can be
taken arbitrarily small while leaving $\Theta_M(k)$ possible $h$.
{Require also that $\epsilon(\delta_i+(s_i+3)\eta)<\eta/2$
for every row $i$.}
For fixed $M$ each active exterior pattern has both selected and
unused margins $\gg_M k$. Indeed its selected count is at least
$h\ge k/(2K)$ for $\epsilon<1/2$, and it is at most
$(k/K)K(\mu_t-\eta)=k(\mu_t-\eta)$, whereas the actual pattern
has $(\mu_t+o(1))k$ positions.

Process the pivot coordinates $i=r,r-1,\ldots,1$ symbolically,
before fixing $h$. Every processed row receives one of three modes.
In a small-total row the eventual nonanchor counts will be $hc^*$.
In a near-full row they will be assigned by lattice matching.
In a deep row the total and its distance below the stock capacity
will exceed an arbitrarily large prescribed fixed multiple of $k$.
Its nonanchor counts will be chosen congruent to $c^*$ modulo $l$.

If an already deep higher row has a nonanchor lower scalar of
unbounded order in coordinate $i$, declare $i$ deep and assign
it to the \emph{lowest} such higher row. Otherwise all higher deep
rows contribute in coordinate $i$ a fixed rational phase $n\psi$
coming from their prescribed residues. At $M$ their contribution is
$M\psi$. If an active exterior pattern is starred in coordinate
$i$, declare it deep and assign an exterior bank. In the remaining
cases all active exterior entries there are named.

Let $P$ be the union of already near-full higher row blocks and
$E$ the union of small-total higher blocks. If $b_{i,n}|_{W_P}$ has
a wild witness, declare $i$ deep and assign that witness. Otherwise
take the fixed form $V$ above. Write $y_{i,n}$ for the sum of the
exterior coordinates with unscaled counts $x_t$. In $C_{Dn}$ use
{the named form $S=(S_{i,n})_n$ given by}
$$
 S_{i,n}
   =D\widetilde y_{i,n}
    +D\sum_{a\in E} c_a^*\widetilde b_{a,i,n}
    +\sum_{a\in P}(D c_a^*-M t_a)\widetilde b_{a,i,n}.
$$
Tildes denote the sum of fixed standard representatives of the
individual labels with the displayed coefficients. Write
$\theta_i=\sum_{a\in P} t_a\widetilde b_{a,i,n}\bmod D$; it was
stabilized independently of $M$. {For a candidate multiplier $h$ in
\eqref{mixmult}, take any $\kappa_P\in qW_P$ and form the candidate
nonanchor counts $hc_P^*+(n-hM)\rho_P+\kappa_P$ on $P$. Use
$hc_E^*$ on $E$ and the prescribed residues on the higher deep rows.
Let $Y_i\in C_n$ be the resulting contribution in coordinate $i$ from
the exterior and all higher stocks, and let $Y_i^*\in C_M$ be the
corresponding contribution in the proposed auxiliary zero. For every
such $\kappa_P$ we have}
\begin{equation}\label{mixres}
\begin{split}
 D Y_i&=hS_{i,n}+n(\theta_i+D\psi)
                     +(D/q)V(\kappa_P)\pmod{Dn},\\
 D Y_i^*&=S_{i,M}+M(\theta_i+D\psi)\pmod{DM},
 \end{split}
\end{equation}
{Indeed, after multiplication by $D$ the term
$(n-hM)\rho_P$ in the candidate count vector contributes}
$n\sum_{a\in P}t_a\widetilde b_{a,i,n}
 -hM\sum_{a\in P}t_a\widetilde b_{a,i,n}$.
The first part is $n\theta_i$ modulo $Dn$, and the second occurs in
$hS_{i,n}$. {The correction $\kappa_P$ contributes $(D/q)V(\kappa_P)$.}
At exponent $M$, the term with coefficient $-M$ cancels the
$M\theta_i$ term modulo $DM$. The rational phase
terms denote integral residues when combined as displayed.
If $S$ has no bounded multiple, declare $i$ deep by progression
avoidance. If $V|_{W_P}\ne0$, also declare $i$ deep, using \eqref{mixdrift}
when $S$ has a bounded multiple. If $V=0$ and $S$ is tame,
\eqref{mixres} gives the exact test $Y_i=ha_i+n\xi_i\pmod n$,
$Y_i^*=a_i+M\xi_i\pmod M$, with $a_i=A/(lD)$ and the same
rational phase $\xi_i$. Explicitly, if
$S_{i,n}=(A+Dn\phi)/l\bmod Dn$, multiplying by $h\equiv1\bmod l$
preserves the phase coefficient $\phi\bmod l$. Division of the
congruence by $D$ gives $a_i=A/(lD)$ and
$\xi_i=\phi/l+\theta_i/D+\psi\bmod1$, with the same expressions at
$M$. For fixed $M$, $ha_i=O_M(k)=o(n)$. A nonzero phase therefore
puts the residue at distance $\gg_M n$ from zero, giving a deep row. With
zero phase $a_i$ is integral. For $a_i\le0$ declare $i$ small-total;
for $a_i>0$ declare it near-full. In the latter case $a_i>d_i^*$
because $T_i^*=[-a_i]_M\le M-1-d_i^*$. This completes the
classification. A designated pure row needs no nonanchor choices
and may be processed by the same rules.

\medskip\noindent\textit{Producing the classified pivot totals.}\quad
Let $P_f$ be the union of all near-full blocks. Matching in this
space means taking
\begin{equation}\label{mixround}
Z_0=f_{P_f}-h f'_{P_f}-(n-hM)\rho_{P_f},\qquad
 \Delta\in Z_0+qW_{P_f},\qquad
             |\Delta_a-\eta k|<\eta k/(4(|I|+1)).
\end{equation}
To see existence, use any basis matrix $B_f$ of $\mathcal B_{P_f}$.
By the capacity identities,
$$
 B_f Z_0=\lambda_f(k-hK)+hK\eta B_f\mathbf1+o(k)+O(h),
 \qquad \lambda_f=\lim B_f(f_{P_f}-n\beta_{P_f})/k .
$$
For clarity let $R_f$ be a right inverse of the independent rows
of $B_f$ (the zero map if there are no rows). Subtracting
$\eta kB_f\mathbf1$ from this display and dividing by $k$ gives
an error bounded in norm by
$\epsilon\|\lambda_f-\eta B_f\mathbf1\|+O_f(1/K)+o(1)$.
Consequently $\eta k\mathbf1+R_fB_f(Z_0-\eta k\mathbf1)$ belongs
to $Z_0+\ker_{\mathbb R}B_f$ and differs from $\eta k\mathbf1$
by at most $\eta k/(8(|I|+1))$ in sup norm. Choose $M$ sufficiently
large first, taking the maximum over the finitely many possible
$P_f$, and then $\epsilon$ sufficiently small. This order makes the
choice of sufficiently late {auxiliary} exponents uniform over their
possible selections. The bounded covering radius of $qW_{P_f}$
in this kernel absorbs the remaining margin and gives \eqref{mixround}.
$Z_0$ is integral because $n\equiv M\pmod D$ and $h\equiv1\pmod D$.
The nonanchor selections in these rows are
$$
 z=h c_{P_f}^*+(n-hM)\rho_{P_f}+\kappa,qquad
 \kappa=Z_0-\Delta\in qW_{P_f}.
$$
They omit $h(f'_{P_f}-c_{P_f}^*)+\Delta$ sites and have $\omega(k)$
available sites on the selected side. Every component of the
omission vector is at least $\eta k/2$.
For each previously used prefix $P$, projection gives
$\kappa_P\in qW_P$. Substitution of the displayed $z_P$ and
$hc_E^*$ gives \eqref{mixres} for any candidate deep counts with their
prescribed residues. This calculation is valid before choosing the
remaining anchors and the individual deep counts.

If a progression was invoked in \eqref{mixres}, discard all $h$ putting the
residual within the prescribed multiple of $k$ of zero. More
explicitly, write
$$
 V(\kappa_P)=V(f_P-n\rho_P)+hV(M\rho_P-f'_P)-V(\Delta_P).
$$
The last two terms are $O_M(k)$ independently of rounding. Translate
the tested interval by the first term and enlarge it by $O_M(k)$
to cover both errors and all forthcoming kernel moves of $O(k)$.
Only $o(k)$ values are discarded. For a tame $S$ with $V|_{W_P}\ne0$,
the macroscopic phase of $Y_i$ equals a rational number plus
$V(\beta_P-\rho_P)/q$. If $V\beta_P$ is irrational the distance
from zero modulo one is positive. For rational $V\beta_P$,
either that phase is nonzero or the additional term
$V(f_P-n\beta_P)/q$ has magnitude with ratio to $k$ tending to
infinity by \eqref{mixdrift}, and is $o(n)$. All errors are $O_M(k)$. These
rows are therefore deep as declared. For $V=0$, kernel moves do
not affect the exact test at all. A finite union of progression
exceptions leaves an $h$ in \eqref{mixmult}. Fix such an $h$ and \eqref{mixround}.

{
In a small-total row the actual total is $h|a_i|$, at least
$hT_i^*$, hence assign the nonanchors exactly $hc^*$ and fill with
anchors. The needed counts are $O_M(k)$, each supply $\omega(k)$.
For a near-full row the actual total is $n-ha_i$ and the required
omission total is $ha_i-d_i-1$. Since $a_i\ge M-T_i^*$,
$$
 ha_i-d_i-1-h(N_i^*-T_i^*)\ge (s_i+2)\eta k
$$
eventually. {Indeed, the left side is at least $h(1+d_i^*)-d_i-1$.
By \eqref{mixmult} and $d_i/k\to\delta_i$, this is at least
$[(s_i+3)\eta-\epsilon(\delta_i+(s_i+3)\eta)+o(1)]k$.
Our choice of $\epsilon$ gives the displayed inequality.}
Thus the additional nonanchor omissions in \eqref{mixround} leave
at least $\eta k$ extra anchor omissions. All omitted counts are
$O_M(k)$ and all starting supplies are $\omega(k)$. This verifies the
anchor slack used by the witness banks. The resulting row total
does not depend on any witness bank, {because $V$ vanishes on the permitted lattice directions in the exact endpoint equation.} The other deep rows have only fixed torsion phases
on all small-total and near-full pivots.
\par}

For each target assigned a wild prefix witness $w_i\in W_P$, lift
a positive multiple to $w'_i\in W_{P_f}$ by prefix compatibility.
The newly added blocks are lower than or equal to the target pivot.
Their lower offsets vanish in that coordinate, even for its own
block. Thus $b_{i,n}(w'_i)$ is a fixed positive multiple of the
same wild scalar. Use $\lfloor\gamma k\rfloor$ independent count
changes $q w'_i$ per target. They adjust all compensating anchors
within each row so its eventual total is fixed. The individual
nonanchors have $\eta k/2$ sites of room on both sides by \eqref{mixround},
and each near-full anchor has at least $\eta k$ of room {by the preceding omission bound}. A fixed sufficiently small $\gamma>0$ makes every
combination legal. In each assigned coordinate the powered
increments have norm tending to infinity, so the switch estimate
produces $r+1$ options separated by the required fixed multiple of
$k$. These moves are in $qW_{P_f}$. {They preserve all exact endpoint equations}, whereas the other nonwild bounded remainders move by
$O(k)$, already allowed for above. {Simultaneously use the
\hyperref[construction:exterior-exchanges]{exterior-pattern exchanges}
for the starred targets assigned to them. The selected and unused margins
verified after \eqref{mixmult} meet their hypotheses for fixed $M$.} {Each assigned coordinate has one bank, either an exterior bank or a witness bank. There are at most $r$ such coordinates,} and each bank's conditional
hit probability is at most $1/(r+1)$ whatever the other options
are. At every target secured in this step, all higher deep rows
contribute only their prescribed phases. Targets assigned to a
higher deep row will be settled when that row is processed.

It remains to realize the deep rows, in descending order. If a target
$i$ was assigned to a particular deep row $j>i$, all subsequent deep
rows above $i$ have torsion-only contrasts into $i$ by the choice of
the lowest such $j$. Their contributions there are fixed in advance
by the prescribed residues. Near-full and small-total lower-offset
contributions there have already been fixed. Thus when processing
$j$, one can calculate in advance the residual at every target
assigned to $j$, apart from the choices in row $j$ itself.
Row $j$ already has a determined total with both margins larger than
an arbitrarily prescribed fixed multiple of $k$. Start with all
nonanchor counts congruent to $c^*$ modulo $l$ and reserve many
disjoint selected-anchor/omitted-nonanchor swaps in groups of $l$.
Explicitly, for each assigned target choose one nonanchor contrast
of unbounded order into it. Reserve a desired fixed multiple of $k$
swaps with that type, while leaving at least $l|I|$ additional sites
selected and omitted in every label. Such reservations are possible
since all labels have $\omega(k)$ sites and both total margins can
exceed twice the sum of the reservation sizes. Fill arbitrarily to
the required total, then adjust each nonanchor count by at most
$l-1$ using the extra anchor room to obtain the prescribed residue.

For a required margin $H_i k$ at target $i$, avoid the circular
arc $[-\lceil H_i k\rceil-1,d_i+\lceil H_i k\rceil+1]$.
Indeed outside this arc the least residue of $-Y_i$ and its
distance below $n-1-d_i$ both exceed $H_i k$. Its number of
residues is at most $(2H_i+4)k$ for large $n$. Request separation
strictly greater than this entire length in every bank assigned
to $i$. For any such fixed separation scale, enough groups give $r+1$
separated choices in each target even if the powered $l$-contrast
has centered magnitude just one: unbounded order implies this fixed
multiple is nonzero eventually, and the sum of capped norms exceeds
the switch bound by using a sufficiently large multiple of $k$
groups. Independent choices put every assigned target at distance greater than the prescribed margin $H_i k$
from both endpoints with positive probability.
Choose the required depth constants from the bottom upward before
the construction: a row needs enough room to supply the finite
separation budgets of its strictly lower assigned targets. The
finite acyclic order gives finite constants throughout. All
unassigned earlier guarantees persist, since if a later deep
source could affect an earlier-guaranteed target by more than its
prescribed phase it would have been that target's assigned source
or one preceding it. This completes the descending {choice of positions in the original sequence}. Each pivot total cancels its residual exactly. The
nonempty exterior selection gives {a zero in the original sequence}, a contradiction.

\section{Realization on indexed positions and proof of Theorem~\ref{thm:main}}\label{sec:positions}
{A nonempty zero sum at a fixed auxiliary exponent $M$ would give a
nonempty zero sum on distinct original positions at all sufficiently
late exponents $n$, even when several coordinates vary within a class.}
We then apply ordinary flag stability to the auxiliary sequences and
recover a flag with smaller deficit in the original sequences
{by comparing the two flags {at the same exponent}}.

\subsection{Realizing the class counts}\label{sec:realize}

{
\medskip\noindent\phantomsection\label{construction:stock-variants}\textit{Auxiliary stock variants.}\quad
Use the dispersed coordinate sets $W_A$ of Subsection~\ref{sec:exterior-exchanges}.
The constant coordinates
of a distinguished class in each row vanish below the row pivot.
In {the auxiliary sequences} use the constants at exponent $M$ and put
zero at the dispersed coordinates to define the base label. For each $i\in W_A$, add the variant obtained by adding $e_i$ to this base label.
The distinguished base label is the original basis vector.
The number of variants is fixed. Each can be given more than
$5K$ positions if its class capacity exceeds $6(r+1)K$, where
$K=M/C$. All computations of how many sites are needed from a
class use the sum across its variants.

To apply the joint count formulas, take $I_j$ to index the
nondistinguished \emph{classes} of row $j$. Set each base lower entry
to its common value if constant and to zero if dispersed. The
distinguished base has lower part zero. Choose the margins in
$C,v$ greater than $10(r+1)$ in place of $10$. This is possible by
the same $o(1)$ correction before fixing $C$. Divide each {auxiliary} class's
stock as evenly as possible among its at most $r$ variants, including
the base. For all large $M$ every variant then has more than $5K$
positions. Use aggregated counts $c_a^*$ in the formulas of
Subsection~\ref{sec:counts}.
\par}

\begin{theorem}\label{aux}
For a critical zero-sumfree family as in Section~\ref{sec:tools}, perform
the preparation in {Lemma~\ref{lem:prepare-stock}} and use its enlarged exterior
size $k$. {After passage to a subsequence, the class data and a constant $C$
can be chosen independently of the finite exterior pattern set $F$ and
the admissible margin $\eta$. For each fixed pair $(F,\eta)$, there is
a further subsequence with the following property.}
With the class capacities of Subsection~\ref{sec:counts} and the
{\hyperref[construction:stock-variants]{stock variants defined above}}, {for every sufficiently large $M$ on this subsequence, any
nonempty auxiliary zero would give a nonempty zero in the original
sequences at all sufficiently late $n$ on it. The threshold on $M$ is
independent of the proposed zero. The threshold on $n$ may depend on
$M$ and on that zero. Hence} the auxiliary sequences are zero-sumfree
{for all sufficiently large $M$}.
\end{theorem}

We give the additional lifting argument needed when a stock class in a
triangular flag has lower entries which vary from position to position.
It identifies the exact requirements on the finite-class count argument.
Write $\|a\|_n$ for circular distance to zero. Throughout, $r$ and the
number of stock classes are fixed, $k\to\infty$, and $k=o(n)$.
The ordinary row $j$ has coordinate $j$ equal to one and all higher
coordinates zero. Its retained population is $n-1-d_j$, $d_j=O(k)$.

\begin{proof}[{Proof of Theorem~\ref{aux}}]
\textit{Which classes can contribute a varying coordinate.}\quad
In the scale transfer fix the {auxiliary} exponent $M$ first and write
$K=M/C$. {Use the multiplier $h$ and the tolerance $0<\epsilon<1/2$ of \eqref{mixmult}, so $(1-\epsilon)k/K\le h\le k/K$.} Thus $h\ge k/(2K)$ as the
{original} exponent tends to infinity with $M$ fixed. For a proposed
{auxiliary} zero denote the selected count of class $A$ by $x_A$,
including all variants. The {auxiliary} total {$T'_j=T_j^*$ in row $j$ is the sum of its class counts.} In the count argument there are three
possible prescriptions for the {original} rows.
\begin{enumerate}
\item At the lower endpoint the required total is $h|a|$ for
a fixed integer $a\le0$, with $|a|\ge T'_j$ and
$|a|\equiv T'_j\pmod M$. For nondistinguished classes take
$hx_A$, and assign the balance to the distinguished class.
\item At the upper endpoint every class has $\gg k$ selected
positions and at least $c_M k$ omissions for some $c_M>0$.
Here $\gg k$ need only be a positive constant multiple of $k$.
A change in the joint class-count lattice is allowed within
these bounds.
\item At a row with adjustable total, both total margins exceed
a chosen large constant times $k$. Every class can be given
at least $c_M k$ selections and omissions, with congruences fixed
on nondistinguished counts. The distinguished count absorbs
the remaining total.
\end{enumerate}
Call a class active if its {original} count is positive. At the lower
endpoint the distinguished count is
$h(x_{A_0}+|a|-T'_j)$, and all other counts are $hx_A$.
Thus each active class has at least $h$ selections there.
Every count at this endpoint is $O_M(k)$, whereas every available
class has size $\omega(k)$, so there are also at least $h$
omissions. An inactive class has $x_A=0$, including the
distinguished one since both summands inside the parentheses
are nonnegative. At the other two row types all classes are active.
These facts determine the active classes before individual {original}
positions are chosen. {Every class used in the auxiliary selection is active in the original sequence. The tests below also include a distinguished class that becomes active because its original count contains an extra multiple of $M$.}

\medskip\noindent\textit{Putting the choices in an order that preserves exact equations.}\quad
Process the symbolic pivot tests from the highest row down.
{Assign the target $i$ to the lowest adjustable higher row that can
vary its contribution there through a dispersed entry or an
unbounded-order contrast of constant entries. If no such row exists,
assign $i$ to an active endpoint class dispersed there, or, if none
exists, to an active exterior pattern starred there. In each of these
cases mark row $i$ adjustable. At every remaining target use the
finite-label count tests of Subsection~\ref{sec:counts}.}
The active classes above $i$ have already been determined, so this rule
does not refer to the yet unknown {original} total in row $i$.
{Where a count test is applied, the active exterior entries are
named and all active higher classes have constant entries.}
Inactive classes have zero count in both {sequences}.
Consequently all sums used in a test at this coordinate agree
exactly with the finite-label calculation using base labels.
This applies also to every lower or upper endpoint: bank choices
never change the incoming coordinate of an endpoint row.

For clarity, the remaining realization has the following order.
The finite-label calculation supplies a common multiplier and
joint endpoint counts, with margins that persist under its allowed
lattice changes. {For a target assigned to a higher adjustable
row, settle its interval when that row is processed.} Book in advance the bounded-order contributions of
all later adjustable rows, using their prescribed class residues.
For coordinates with no such higher adjustable row, first settle
the lattice variations in the counts at upper endpoints at those
coordinates where all {active} classes have constant
entries and the exterior is named. A sum at such a coordinate is
an exact function of the counts regardless of the {original} positions
ultimately used. The finite-label progression estimates and the
{exact endpoint equations} apply there unchanged. Exterior switches
also leave these named entries fixed.

With the resulting endpoint counts fixed, choose their {original}
positions and all exterior positions. {At each target assigned to an
endpoint class or an exterior pattern, use a bank in that class or
pattern.} {The active endpoint classes have the
selected and omitted margins established above, and the same is true
for every selected exterior pattern by \eqref{mixmult} and the estimates
following it. Lemma~\ref{pos-banks} and the
\hyperref[construction:exterior-exchanges]{exterior-pattern exchanges}
therefore apply jointly.}
Their sums include any fixed contribution from other classes,
and the booked torsion terms. No bank alters a previously settled
coordinate unless a dispersed class or active exterior star also
exposes it. Those coordinates were reserved for this step or for
an adjustable row.

Finally process the adjustable rows in descending order, with the
incoming required total now known when a row is reached. To avoid
the intervals of targets assigned to this last row, use a
within-class bank at each target with an active dispersed class.
At any other assigned target all entries of the row are constant
within their class. Use the finite-label count exchanges of a
nonanchor class against the anchor, in groups of a fixed size
preserving the prescribed bounded-order residues. Reserve all
these banks on disjoint sites before choosing fillers. A
within-class bank costs $\lfloor\gamma k\rfloor$ chosen and unused
sites for arbitrarily small fixed $\gamma>0$. The interclass
banks cost $O(k)$ with constants fixed by the required target
depths, the rank and the residue modulus. Reserve the sites,
choose a filler count in each class matching its required
congruence, and distribute the remaining total among classes.
For this last assertion each nondistinguished class offers an
arithmetic progression with bounded step, and the distinguished
class offers an integer interval of length $\omega(k)$.
Equivalently fix an anchor residue matching the total; sums of
the resulting arithmetic progressions of the same step cover all
values of that residue between their extreme sums. The discrepancy
at either endpoint is bounded by the number of classes times the
step. The required row margins exceed the bank costs plus that
discrepancy. The construction can also leave at least a fixed
positive multiple of $k$ chosen and omitted in each class by
including those small reserves in the same calculation.

A within-class bank preserves all class counts. Within an
interclass bank toggling a block changes each participating count
by a multiple of the residue modulus. All its increments in its
own target have the prescribed constant value independent of the
particular sites. There are consequently $r+1$ separated options
in the assigned target of each bank (for an interclass bank use
the cyclic switch estimate with sufficiently many groups).
Choose independently and use the union bound of Lemma~\ref{pos-banks}.
The desired depths can be chosen successively from the bottom row
upwards: for a source row require greater depth than the total
bank costs needed to enforce its lower target depths. This is a
finite recurrence independent of $n$. Once an assigned target is
settled, every subsequent adjustable higher row contributes there
only its booked bounded-order offsets. Thus the resulting choices
are simultaneous choices of indexed positions, even when entries
of a class are arbitrarily correlated across coordinates.

At each coordinate where \eqref{mixres}
is invoked, the descending classification ensures that every {active} higher class has a constant entry there. Each inactive higher
class has zero count in both {sequences}. Thus replacing the higher classes
by their base labels preserves both identities in \eqref{mixres} at
each such coordinate. At a target settled by a later adjustable row, the
positions chosen earlier enter the residual with their actual contributions.
Bounded-order offsets at each unexposed target of an adjustable row
are computed using the aggregate class counts modulo $l$. Its
distinguished class has zero constant entry there. The classification of
Subsection~\ref{sec:counts} gives the three count prescriptions, with the
additional adjustable targets assigned as above. In particular, it
requires no estimate on the unknown full sum of a dispersed class;
its target has been assigned an exchange whenever that class is used.
The preceding construction produces {a nonempty zero in the original sequence} from
any auxiliary zero. This proves Theorem~\ref{aux}.
\end{proof}

\subsection{Recovering the flag and proving eventual equality}\label{sec:finish}

We first prove a lemma comparing the rows of two flags.

\begin{lemma}\label{pivots}
Fix $r$ and $C\ge1$, let $M\to\infty$ and put $K=M/C$.
Let a {sequence} of length $rM+o(M)$ in $C_M^r$ contain disjoint
stocks $S_j$, $1\le j\le r$. In a basis $e_1,\ldots,e_r$ each
$S_j$ lies in $e_j+\langle e_1,\ldots,e_{j-1}\rangle$, contains
$e_j$, and has size $M-K\delta_j+o(K)$, where $0\le\delta_j\le1$.
Each distinct type within each stock occurs there more than $5K$ times.
Suppose some basis $f_1,\ldots,f_r$ has $M-o(M)$ or more terms in
each of its ordinary rows $R_h=f_h+\langle f_1,\ldots,f_{h-1}\rangle$.
For large $M$ there is a permutation $j\mapsto h_j$ with
$S_j\subseteq R_{h_j}$ and $e_j\in R_{h_j}$. The two bases $f$ and
the list $e_j$ in increasing $h_j$ order give the same individual
ordinary rows. In particular the mass of terms outside the stocks
lying in $R_{h_j}$ is $K\delta_j+o(K)$.
\end{lemma}
\begin{proof}
The $R_h$ are disjoint. Each contains at most $M+o(M)$ terms,
by subtracting the lower bounds for the other $r-1$ from the
total length. The number outside their union is $o(M)=o(K)$.
Consequently every stock type is in one of the $R_h$ for large $M$.
For $S_1$ the only type is $e_1$, which defines $h_1$.

Suppose the first $j-1$ stocks have been assigned distinct rows.
Each assigned row already contains at least $M-K+o(K)$ stock
terms. Thus it has space for at most $K+o(K)$ further terms, which
excludes each type of $S_j$. All its types lie at unassigned pivots.
In $f$-coordinates, suppress the coordinates at pivots $h_i$, $i<j$,
in descending order by subtracting the corresponding multiple of
$e_i$. Here the pivot coefficient in $e_i$ is one and every higher
coordinate is zero. The reduction is linear. Its kernel is the
span of those $e_i$: it kills them all by this triangular rule, and
the difference between a vector and its reduced vector lies in
their span. All elements of $S_j$ therefore have the same reduced
vector, since they are equal modulo that span. If a vector of $S_j$
belongs to $R_h$, then no occupied pivot above $h$ is used to reduce
it, since all those coordinates are zero. Reductions at lower
pivots do not change its $h$th coordinate or any higher one. Thus
the common reduced vector still has highest nonzero coordinate $h$
and coefficient one there. This proves that all types in $S_j$ lie
at the same previously unassigned pivot $h_j$, including $e_j$.
Induction gives the permutation. In its order the $e_j$ have a
unit triangular change of basis from the $f_h$. Each initial
span agrees, and $e_j-f_{h_j}$ lies in the common preceding span.
This proves equality of individual rows, also over a composite modulus.
Subtracting $|S_j|$ from the full size $M+o(M)$ of this row, which
contains no other stock, gives the asserted count.
\end{proof}

\Needspace{6\baselineskip}Now we are in a position to prove Theorem~\ref{thm:main}.

\begin{proof}[Proof of Theorem~\ref{thm:main}]
{For each prime $p\mid n$, let $v_p(n)$ be the largest
integer $a$ such that $p^a\mid n$.} Write $P(n)=\max_{p\mid n}p^{v_p(n)}$. Suppose Theorem~\ref{thm:main}
fails at a fixed $r$. The basis {sequence} $\prod_{j=1}^r e_j^{\,n-1}$
is zero-sumfree: each of its subsums has its $j$th coordinate among
$0,\ldots,n-1$, and a zero requires every count to vanish. Consequently
the lower bound is $T_r(n)$ for every $n$, and failure gives infinitely
many distinct $n$ with a zero-sumfree critical sequence.
For each such $n$ minimize the exterior size after capping every
ordinary row at $n-1$ over all ordered bases. Denote the minimum by
$k_0$. Since $T_r(n)/n\to r$, Theorem~\ref{stability}
gives $k_0=o(n)$ on a subsequence. In a minimizing basis the sum of
the capped deficits is $k_0-1$. Applying Theorem~\ref{small-main}
with $Q=P(n)$ to this zero-sumfree sequence gives $k_0-1\ge P(n)$.
As $n\to\infty$,
$P(n)\to\infty$: if all prime powers dividing $n$ were at most $B$,
then $n$ would divide $\operatorname{lcm}(1,\ldots,\lfloor B\rfloor)$.
Thus $k_0\to\infty$.

Choose a minimizing basis. {Lemma~\ref{lem:prepare-stock}} moves at most
$O(k_0)$ positions to the exterior, leaving $k_0\le k=O(k_0)$.
In the prepared basis $e_1,\ldots,e_r$, write $n-1-d_j$
for the retained stock in row $j$. Since the length is $r(n-1)+1$,
$k=1+\sum_jd_j$. Pass to a subsequence on which $d_j/k\to\delta_j$.
Thus $0\le\delta_j\le1$ and $\sum_j\delta_j=1$. Take the exterior
pattern masses $\mu_t$ in this same basis as in {Subsection~\ref{sec:exterior-exchanges}}.
They sum to one. Fix the constant $C$ of Theorem~\ref{aux} before
choosing the finite sets below.\par

Exhaust the exterior distribution by nested finite sets $F_s$ of
positive-mass patterns. Choose $\eta_s>0$ smaller than their least
mass with {$(|I|+3)\eta_s<1/4$ at every stage and}
$(|F_s|+|I|+r+3)\eta_s\to0$. At stage $s$ use a nested
subsequence on which all finite tests stabilize.
For this fixed $F_s,\eta_s$, choose {an exponent from the original family}
$M_s>M_{s-1}$ sufficiently late for Theorem~\ref{aux} and the
{\hyperref[construction:stock-variants]{variant multiplicity bounds}}. Also require $|F_s|/M_s<1/s$ and
$\sum_{t\in F_s}|p_s(t)-\mu_t|<1/s$, where $p_s(t)$ is the number
of {original} exterior positions of pattern $t$ at $M_s$ divided by
the {original} exterior size $k$. These requirements concern finitely
many convergent frequencies at each stage. Put $K_s=M_s/C$.\par

By Subsection~\ref{sec:counts}, the auxiliary exterior
count of pattern $t\in F_s$, divided by $K_s$, is
$q_s(t)=\lfloor K_s(\mu_t-\eta_s)\rfloor/K_s$. There are no
auxiliary exterior positions of patterns outside $F_s$. In particular,
$$
 \sum_{t\in F_s}|p_s(t)-q_s(t)|
 \le 1/s+|F_s|(\eta_s+1/K_s)=o(1).
$$
Since $\sum_{t\in F_s}\mu_t\to1$, the auxiliary exterior has
$K_s+o(K_s)$ positions. The auxiliary stock $S_j$ has size
$M_s-1-\lceil K_s(\delta_j+(s_j+3)\eta_s)\rceil
=M_s-K_s\delta_j+o(K_s)$ along this diagonal sequence, where $s_j$
is the fixed number of nondistinguished classes in row $j$. Summing
over $j$ and adding the exterior gives total length $rM_s+o(M_s)$.\par

Theorem~\ref{aux} gives zero-sumfreeness, so
Theorem~\ref{stability} gives, on a further subsequence, a basis with
$M_s+o(M_s)$ positions in each ordinary row. By {the
\hyperref[construction:stock-variants]{construction of the stock variants}}, the auxiliary stocks are disjoint, each
lies in its original $e$-row and contains its distinguished base
$e_j$, and each stock type occurs there more than $5K_s$ times.
Lemma~\ref{pivots} now applies. Write $h_j$ for the new pivot of
stock $j$ and ${H'_j}=\langle e_i:h_i<h_j\rangle$. The recovered row
is $e_j+{H'_j}$ and contains $K_s\delta_j+o(K_s)$ auxiliary exterior
positions. We use this same ordering of the prepared basis in the
{original sequence} at exponent $n=M_s$.\par

First consider {an original} stock class $A$ of row $j$.
Its auxiliary base $b_A$ copies all constant coordinates at $M_s$
and is zero in its dispersed coordinates $W_A$. Both $b_A$ and
$b_A+e_i$ lie in $e_j+{H'_j}$ for every $i\in W_A$, so $e_i\in {H'_j}$.
Every {original} term $x$ of $A$ therefore satisfies
$x-b_A=\sum_{i\in W_A}x_i e_i\in {H'_j}$. Hence all $n-1-d_j$ {original}
stock positions of row $j$ lie in $e_j+{H'_j}$.\par

Next let $A_{j,s}$ consist of the patterns in $F_s$
whose auxiliary vector lies in $e_j+{H'_j}$. {For each $t\in A_{j,s}$, the auxiliary vector encoding $t$ has coefficient one at $e_j$ and coefficient zero at every $e_i$ with $h_i>h_j$.}
A star is encoded by $2$, distinct from zero and one for $M_s>2$.
Thus the pivot and every later coordinate are named, and their
entries agree with those of every {original} occurrence of the pattern
at this same exponent. Each such occurrence belongs to $e_j+{H'_j}$.
Counting indexed exterior positions gives
$\sum_{t\in A_{j,s}}q_s(t)=\delta_j+o(1)$. The preceding absolute
error bound applies to each subset $A_{j,s}$, even when the ordering
varies with $s$. Thus the row contains at least $k(\delta_j-o(1))$
{original} exterior positions.\par

Adding its retained stock, each new {original} row
contains at least $n-1-d_j+k\delta_j-o(k)=n-1-o(k)$ positions.
After capping at $n-1$, each row has deficit $o(k)$. {Let $k_{\mathrm{new}}$ be the capped exterior size of the original
sequence at $n=M_s$ in this reordered basis. Since $r$ is fixed, we have
$k_{\mathrm{new}}=1+o(k)=o(k_0)$, using $k=O(k_0)$ and $k_0\to\infty$.
Thus $k_{\mathrm{new}}<k_0$ for all sufficiently large $s$ along this
subsequence. This contradicts $k_{\mathrm{new}}\ge k_0$, since $k_0$ is
the minimum over all ordered bases for this same original sequence.}
There are
only finitely many exceptional exponents at this fixed rank. Their
maximum, or $1$ when the set is empty, gives the asserted $N(r)$. This completes the proof of the theorem.
\end{proof}

\endgroup

\end{document}